\documentclass[review,1p,times]{elsarticle}

\usepackage{graphicx}

\usepackage{amsmath, amssymb, amsthm, dsfont, bm}         
\usepackage{amsfonts}         
\usepackage{color}

\def\vec#1{\mathbf{#1}}

\newcommand {\R} {{\mathds R}}

\newcommand {\C} {{\mathds C}}

\newcommand{\ie}{{\em i.\thinspace{}e. }}

\newcommand{\eg}{{\em e.\thinspace{}g. }}

\newcommand{\pl}{\partial }

\newcommand{\hb}{\hbar}
\newcommand{\agrad}{ A^t \nabla}
\def\rtf{\rho_{\hbox{\tiny TF}}}
\def\aho{a_{\hbox{\footnotesize ho}}}
\def\omegap{\omega_{\perp}}
\def\vadim#1{\tilde{#1}}
\def\pscal#1#2{{\langle #1, #2 \rangle}}
\def\mod#1{{\left| #1 \right|}}

\def\grdGn{{\mathcal G}_n}
\def\Real{\Re}
\def\Blue#1{{\color{blue}{#1}}}
\def\Red#1{{\color{red}{#1}}}

\journal{Computer Physics Communications}

\begin{document}

\begin{frontmatter}



\title{A finite-element toolbox for the stationary Gross-Pitaevskii equation with rotation}


\author[lmrs,ljll]{Guillaume Vergez}
\ead{vergez@ann.jussieu.fr}

\author[lmrs]{Ionut Danaila\corref{io}}
\ead{ionut.danaila@univ-rouen.fr}

\author[ljll]{Sylvain Auliac}
\ead{auliac@ann.jussieu.fr}

\author[ljll]{Fr{\'e}d{\'e}ric Hecht}
\ead{hecht@ann.jussieu.fr}

\address[lmrs]{Universit{\'e} de Rouen, Laboratoire de Math{\'e}matiques Rapha{\"e}l Salem,  CNRS UMR 6085, Avenue de l'Universit{\'e}, BP 12, F-76801 Saint-{\'E}tienne-du-Rouvray, France}

\address[ljll]{Sorbonne Universit{\'e}s, UPMC Univ Paris 06, CNRS UMR 7598, Laboratoire Jacques-Louis Lions, 4 Place Jussieu, F-75005 Paris, France}

\cortext[io]{Corresponding author. Tel.: (+33)2 32 95 52 50 ;
              fax: (+33)2 32 95 52 86}

\begin{abstract}
We present a new numerical system using classical finite elements with mesh adaptivity for computing stationary solutions of the Gross-Pitaevskii equation. The programs are written  as a toolbox for FreeFem++ (www.freefem.org), a free finite-element software available for all existing operating systems. This offers the advantage to hide all technical issues related to the implementation of the finite element method, allowing to easily implement various numerical algorithms.Two robust and optimised numerical methods were implemented to minimize the Gross-Pitaevskii energy: a steepest descent method based on Sobolev gradients and a minimization algorithm based on the state-of-the-art optimization library Ipopt. For both methods, mesh adaptivity strategies are implemented to reduce the computational time and increase the local spatial accuracy when vortices are present. Different run cases are made available for 2D and 3D configurations of Bose-Einstein condensates in rotation. An optional  graphical user interface is also provided, allowing to easily run predefined cases or with user-defined parameter files.  We also provide several post-processing tools (like the identification of quantized vortices) that could help in extracting physical features from the simulations. The toolbox is extremely versatile and can be easily adapted to deal with different physical models. 
\end{abstract}

\begin{keyword}
FreeFem++  \sep Ipopt \sep Gross-Pitaevskii \sep Bose-Einstein   \sep finite element \sep mesh adaptivity \sep Sobolev gradient.

\end{keyword}

\end{frontmatter}
   \noindent
   {\bf Programm summary}\\
  {\em Program Title:} GPFEM\\
{\em Catalogue identifier:}\\ 
{\em Program summary URL:}\\
{\em Program obtainable from:}\\
{\em Licensing provisions:} Standard CPC licence, http://cpc.cs.qub.ac.uk/licence/licence.html\\
{\em No. of lines in distributed program, including test data, etc.:} 49177\\
{\em No. of bytes in distributed program, including test data, etc.:} 450969\\
{\em Distribution format:} tar.gz\\
{\em Programming language:} FreeFem++ (free software, www.freefem.org)\\
{\em Computer:} PC, Mac, Super-computer.\\
{\em Operating system:} Windows, Mac OS, Linux.\\
{\em Classification:} 2.7, 4.9, 7.7.\\
{\em Nature of problem:} The software computes 2D or 3D  stationary solutions of the Gross-Pitaevskii equation with rotation. The main application is the computation of different types of vortex states (Abrikosov vortex lattice, giant vortex) in rotating Bose-Einstein condensates. The software can be easily modified to take into account different related physical models. \\
{\em Solution method:} The user has the choice between two robust and optimised numerical methods for the direct minimization of the Gross-Pitaevskii energy: a steepest descent method based on Sobolev gradients and a minimization algorithm based on the state-of-the-art optimization library Ipopt. For both methods, mesh adaptivity strategies are implemented to reduce the computational time and increase the local spatial accuracy when vortices are present. \\
{\em Running time:}\\ From minutes for 2D configurations to hours for 3D cases (on a personal laptop).  Complex 3D cases (with hundreds of vortices) may require several days of computational time.\\


\pagebreak

\section{Introduction}

The Bose-Einstein condensate (BEC) is an ideal system to study superfluidity at a macroscopic level: it is a highly controllable quantum system which admits a simple theoretical description using the Gross-Pitaevskii  equation (GPE) \cite{BEC-book-2003-pita}. A great deal of attention has been lately devoted to the development of accurate numerical schemes to solve different forms of the GPE, from the classical (stationary or time-dependent) GPE,  to systems of coupled GPEs (\eg for two-component or spinor BEC) and more recent formulations (\eg with non-local interactions or fractional GPE). For recent reviews of numerical methods for GPE, see 
\cite{BEC-review-2004-Minguzzi,BEC-review-2006-bao,BEC-review-2013-bao-KRM,BEC-review-2013-antoine-besse-bao,BEC-review-2014-bao-ICM}.

Among all these formulations, the stationary GPE is used either to numerically generate an initial condition for the simulation of real-time dynamics of BEC, or to directly investigate physical features of experimentally observed BEC. In the former case, the stationary (ground state) solution which is the global minimizer of the GP energy is sought, while in the latter case, capturing local minima of the GP energy could be of interest since they represent  excited (or metastable)  states observed in experimental BEC configurations. The most striking example of how numerical solutions of the stationary GP equation were used to investigate physics is the study of quantized vortices in rotating BEC. Since superfluidity in BEC is closely related to the nucleation of quantized vortices, this topic has focused the attention of physical and mathematical communities during the last two decades. Numerous experimental and theoretical studies were devoted to the investigation of three-dimensional properties of single (straight or bent) vortex lines, vortex rings or Abrikosov vortex lattices (for a review of such physical systems, see the dedicated volumes
\cite{BEC-book-2001-barenghi,BEC-book-2008-panos,BEC-book-2008-barenghi,BEC-book-2009-tsubota}). Numerical simulations of the stationary three-dimensional (3D) GPE proved as a valuable investigation tool for all these topics, revealing properties of quantized vortices difficult to observe experimentally, suggesting new configurations, or supporting new physical or mathematical theories (\eg \cite{BEC-numm-2005-kasamatsu-3d,BEC-numm-2004-berloff-3d,dan-2003-aft,dan-2004-aft,dan-2005}; for a review, see \cite{BEC-review-2008-kasamatsu}). 


The difficulty in computing solutions of the stationary GP equation with rotation comes from the presence in a condensate of a large number of vortices, with large gradients of atomic density in the vortex cores. This explains the use in the literature of discretisation methods with high order spatial accuracy: Fourier spectral \cite{BEC-physV-2001-perez-gc,BEC-physV-2001-perez-gcM,BEC-CPCm-2009-zeng}, sixth-order finite differences \cite{dan-2003-aft,dan-2004-aft,dan-2005}, sine-spectral \cite{BEC-baow-2004-Du,BEC-baow-2006-gse}, Laguerre--Hermite pseudo-spectral \cite{BEC-baow-2008-sym}, hybrid discontinuous Galerkin discretisations based on polynomials and plane waves \cite{BEC-numm-2012-Farhat}, etc. Several software packages for the classical stationary GPE  
were deposited in the CPC Program Library. They use different numerical methods: iterative diagonalization method  \cite{BEC-CPCs-2006-Tiwari}, optimal damping algorithm \cite{BEC-CPCs-2007-dion-cances,BEC-CPCs-2014-Hohenester}, Crank-Nicolson scheme \cite{BEC-CPCs-2009-Muruganandam}, Newton-like method with an approximate line-search
strategy \cite{BEC-CPCs-2013-Caliari,BEC-CPCs-2012-Vudragovic}, fully-explicit fourth-order Runge–Kutta scheme \cite{BEC-CPCs-2013-Caplan}, semi-implicit backward Euler scheme  \cite{BEC-CPCs-2014-antoine-duboscq}, etc. The spatial discretization is generally based on spectral \cite{BEC-CPCs-2007-dion-cances,BEC-CPCs-2013-Caliari,BEC-CPCs-2014-antoine-duboscq} or finite-difference \cite{BEC-CPCs-2009-Muruganandam,BEC-CPCs-2012-Vudragovic,BEC-CPCs-2013-Caplan,BEC-CPCs-2014-Hohenester} methods. Provided programs are written in Fortran  \cite{BEC-CPCs-2007-dion-cances,BEC-CPCs-2009-Muruganandam}, C \cite{BEC-CPCs-2012-Vudragovic,BEC-CPCs-2013-Caplan} or Matlab \cite{BEC-CPCs-2013-Caliari,BEC-CPCs-2013-Caplan,BEC-CPCs-2014-antoine-duboscq,BEC-CPCs-2014-Hohenester}.

Numerical methods based on standard finite elements are less represented in this field. Vortex states in  rotating BEC were computed using finite elements with fixed meshes \cite{BEC-math-2001-Aftalion-Du,BEC-baow-2003-gs,BEC-baow-2004-Du,BEC-numm-2009-baksmaty} or dynamically adapted meshes \cite{dan-2010-JCP}, but only for 2D configurations.  To the best of our knowledge, no finite-element programs exist in the CPC Program Library for the GP equation with rotation. The purpose of this paper is thus to distribute a finite-element solver for computing steady solutions of the GPE with rotation, in both 2D and 3D settings. The code was built as a toolbox for FreeFem++ \cite{hecht-2012-JNM,freefem}, which is a free software (under LGPL license) using a large variety of triangular finite elements  (linear and quadratic Lagrangian elements, discontinuous $P^1$, Raviart-Thomas elements, etc.)  to solve partial differential equations. FreeFem++ is an integrated product with its own high level programming language and a syntax close to mathematical formulations, making the implementation of numerical algorithms very easy. Among the features making FreeFem++ an easy-to-use and highly adaptive  software we recall the advanced automatic mesh generator, mesh adaptation, problem description by its variational formulation, automatic interpolation of data, colour display on line, postscript printouts, etc. FreeFem++ community is continuously growing, with  thousands of users all over the world.

The present FreeFem toolbox, called GPFEM, provides two efficient numerical methods for computing stationary states with vortices, with the following novelties:\\
{\em (i)}  the steepest-descent algorithm based on Sobolev gradients suggested in \cite{dan-2010-SISC} and tested for 2D configurations in \cite{dan-2010-JCP} was improved by adding an optimized line-search algorithm for the descent step and extended for 3D configurations;\\
{\em (ii)} a novel minimisation method for 2D and 3D configurations was implemented based on  the state-of-the-art optimisation library Ipopt \cite{Ipopt-int-point-2002} using the direct minimization interior point method;\\
{\em (iii)} the mesh adaptivity algorithm suggested in \cite{dan-2010-JCP} for 2D configurations was extended in 3D  and optimised by the use of anisotropic mesh adaptivity functions provided by \textbf{mshmet} \cite{mshmet} and \textbf{mmg3d} \cite{mmg3d} softwares.\\
From the programming point of view, the toolbox presents the following advantages:\\
{\em (iv)} the switch from different finite elements (from linear $P^1$ to quadratic $P^2$ and high-order $P^3$ or $P^4$ finite elements) implies the modification of a single instruction (the definition of the finite-element space);\\
{\em (v)}  the scripts are easy to adapt for different mathematical  or physical settings  (two different scalings are implemented);\\
{\em (vi)}  a graphical interface allows to run predefined 2D or 3D examples.

The paper is organised as follows. In \S2, we present different  mathematical formulations of the GP equation and energy.  Two different scalings are introduced. Numerical methods are presented in \S3 and the important issue of setting the initial condition for the computation is  described in \S4. The details of the derivation of closed formulae for the Thomas-Fermi approximation (generally used as initial condition) is deferred to  Appendix A. The structure of the provided software is described in great detail in \S5. Various test cases for computing 2D and 3D configurations with vortices are presented in \S6. The optional user interface is also described in \S6. The main features of the software and possible extensions are summarised in \S7. 

\pagebreak

\section{Mathematical model: the Gross-Pitaevskii energy}\label{section-math}

\subsection{The Gross-Pitaevskii energy for the rotating condensate}

We consider in this paper numerical methods for the direct minimization of the Gross-Pitaevskii energy. For a pure BEC of $N$ atoms confined in a trapping potential $V_\text{trap}({\vec x})$ rotating with angular velocity ${\vec \Omega}$, the energy of the system in the rotating 
frame is described by the  functional:
\begin{equation}
\label{eq-GP-energ-gen}
{\mathcal E}(\psi) = \int_{\R^3} \left(\frac{\hb^2}{2m} |\nabla \psi|^2 + V_\text{trap}\, |\psi|^2 +{\frac{1}{2} g} |\psi|^4 -{{\psi^* {\vec \Omega}}}{\cdot} {\vec {\mathcal L}}(\psi)\right)\, d{\vec x}
\end{equation}
where $\psi({\vec x})$ is the classical field complex wave function, $\psi^*$ denotes its complex conjugate, $m$ is the atomic mass, $\hbar$ the reduced Planck constant  and  $g$ the coupling constant
\begin{equation}
\label{eq-GP-g}
g = 4\pi\hb^2a_s/m, \quad \mbox{with $a_s$ the scattering length.}
\end{equation}
The kinetic momentum $\vec {\mathcal L}$ can be expressed as 
\begin{equation}
\label{eq-GP-L}
{\vec {\mathcal L}}(\psi) = {\vec x} \times {\vec P}(\psi), \quad \mbox{with the impulse}\quad {\vec P}(\psi) = -i \hb \nabla \psi,
\end{equation}
We consider in the following rotations along the $z$-axis (\ie ${\vec \Omega}= \Omega\, {\vec k} )$  and therefore only the  $z$-component of the kinetic momentum appears in (\ref{eq-GP-energ-gen}) for the rotation term:
\begin{equation}
\label{eq-GP-Lz}
{\vec \Omega} \cdot {\vec {\mathcal L}}(\psi) = \Omega {\mathcal L}_z \psi =  i \hb \Omega \,\left(y\frac{\pl \psi}{\pl x} - x\frac{\pl \psi}{\pl y} \right) = i \hb \Omega\, \agrad \psi, \quad \mbox{with}\quad
A^t = (y , -x , 0).
\end{equation}
As a consequence, the form of the Gross-Pitaevskii energy considered in this paper is:
\begin{equation}
\label{eq-GP-energ}
{\mathcal E}(\psi) = \int_{\R^3} \left(\frac{\hb^2}{2m} |\nabla \psi|^2 + V_\text{trap}\, |\psi|^2 +{\frac{1}{2} g} |\psi|^4 \right)\, d{\vec x} - \Omega L_z,
\end{equation}
with $L_z$ the total angular momentum:
\begin{equation}
\label{eq-GP-Lzint}
L_z = i \hb \int_{\R^3} \psi^* \agrad\psi d{\vec x}= i \frac{\hb}{2} \int_{\R^3} \left(\psi^* \agrad\psi -
\psi\agrad\psi^* \right) d{\vec x}= \hb \int_{\R^3} \Real\left(i \psi^* \agrad \psi \right) d{\vec x}.
\end{equation}
$\Real$ denotes the real part. We compute here minimizers of the energy (\ref{eq-GP-energ}) with the constraint 
\begin{equation}
\label{eq-GP-cons}
\int_{\R^3} \left|\psi({\vec x})\right|^2\, d{\vec x} = N,
\end{equation}
expressing the conservation of the number of atoms in the condensate. Among these minimizers, the ground state is defined as a global minimum, \ie $\psi_g =\min_\psi {\cal E(\psi)}$. Local minimizers with energy larger than that of the ground state are called excited states or meta-stable states.

Using (\ref{eq-GP-Lzint}), the energy (\ref{eq-GP-energ}) can be written in the following form that will be useful in deriving numerical methods in the next section:
\begin{equation}
\label{eq-GP-energ-Agrad}
{\mathcal E}(\psi) = \int_{\R^3} \left(\frac{\hb^2}{2m} 
\left|\nabla \psi + i \frac{m \Omega}{\hbar}\, A^t \psi \right|^2 + V^\text{eff}_\text{trap}\, |\psi|^2 +{\frac{1}{2} g} |\psi|^4 \right)\, d{\vec x},
\end{equation} 
where the effective trapping potential is corrected with the centrifugal term:
\begin{equation}
\label{eq-GP-Veff}
V^\text{eff}_\text{trap} = V_\text{trap} - \frac{1}{2} m \Omega^2 r^2, \quad r^2 =x^2 + y^2.
\end{equation}

Another useful form of the energy corresponds to the grand potential of the system:
\begin{equation}
\label{eq-GP-energ-grandp}
\Xi = {\cal E}(\psi) - \mu N = {\cal E}(\psi) - \mu \int_{\R^3} |\psi|^2 \, d{\vec x},
\end{equation}
where $\mu \in \R$ is the chemical potential of the condensate, introduced as a Lagrange multiplier 
for the constraint (\ref{eq-GP-cons}). The  Euler-Lagrange equation ($\delta \Xi=0$) corresponding to (\ref{eq-GP-energ-grandp}) leads to the stationary (or time-independent) GP equation:
\begin{equation}
\label{eq-GP-stat}
-\frac{\hb^2}{2m} \nabla^2 \psi + V_\text{trap} \psi + g |\psi|^2 \psi - i \hb \Omega \agrad \psi = \mu \psi.
\end{equation}
The ground state and excited states are therefore eigenfunctions of the nonlinear eigenvalue problem (\ref{eq-GP-stat}).

We also consider in this paper two-dimensional (2D) configurations corresponding to disk-shape (or pancake) condensates. The dimension reduction from 3D to 2D can be done by approximating the 3D wave function by a factorized ansatz $\psi(x,y,z) = \psi_{2D}(x,y) \psi_3(z)$. For the  precise form of the ansatz, the reader is referred to review papers \cite{BEC-review-1999-dalfovo,BEC-review-2015-Bagnato}. For  a mathematical justification of the dimension reduction from 3D to 2D equations, see  \cite{BEC-review-2013-bao-KRM}. By integrating out the $z$-dependence, previous forms of energy and stationary GPE stand, with $\R^3$ replaced by $\R^2$, with the caution that the non-linear interaction constant $g$  expressed by (\ref{eq-GP-g}) for the 3D setting has to be replaced by its {\em reduced} form in 2D. This constant will be prescribed as an input parameter of the computation.

\subsection{Scaling and trapping potential}\label{sec-scal}
We consider in the following the Gross-Pitevskii model set on $\R^d$, with $d=3$ or 2. Various forms of scaling are used in the literature \cite{BEC-physV-2005-fetter,BEC-physV-2002-kasamatsu-2d,BEC-math-2001-Aftalion-Riviere}. To allow the switch between different scalings, we introduce a parameter $\varepsilon$ and define a general length scale as:
\begin{equation}
\label{eq-scal-xs}
\quad x_s = \frac{\aho}{\sqrt{\varepsilon}}, \quad \aho=\sqrt{\frac{\hb}{m \omegap}},
\end{equation}
where $\aho$ is the harmonic oscillator length defined with respect to a reference trapping frequency $\omegap$. By setting $\vadim{\vec{x}} = {\vec x}/{x_s}$ and 
\begin{equation}
\label{eq-scal-u}
u = \frac{\psi}{\sqrt{N} \, x_s^{-d/2}} = \varepsilon^{-d/4}\, \frac{\psi}{\sqrt{N} \, \aho^{-d/2}},
\end{equation}
the dimensionless GP energy (per particle) becomes:
\begin{equation}
\label{eq-scal-energ}
\displaystyle {E}(u) = \frac{{\mathcal E}(\psi)}{N \frac{\hb^2}{m} \aho^{-2}} =
\frac{{\mathcal E}(\psi)}{ N  \hbar \omega_\perp} =
\varepsilon \, \int_{\R^d} \left[ \frac{1}{2} |\nabla u|^2 + C_\text{trap}\, |u|^2 + \frac{1}{2} C_g |u|^4 - i C_\Omega\, u^* \agrad u \right] \, d\vadim{\vec{x}},
\end{equation}
where
\begin{eqnarray}
\label{eq-scal-Ctrap}
\displaystyle C_\text{trap}(\vadim{x},\vadim{y},\vadim{z}) =  \frac{1}{\varepsilon^2} \vadim{V}(\vadim{x},\vadim{y},\vadim{z}), & \displaystyle  \vadim{V}(\vadim{x},\vadim{y},\vadim{z})= \frac{1}{m \omegap^2 x_s^2}\, V_\text{trap}(x,y,z)\\
\label{eq-scal-Cg}
C_g =  \sqrt{\varepsilon} \, \beta,  & \displaystyle   \displaystyle \beta = \frac{4\pi N a_s}{\aho} \quad \text{(in 3D)}, \quad \beta=\beta_{2D}\quad \text{(given in 2D)},\\
\label{eq-scal-Com}
C_\Omega = \displaystyle \frac{1}{\varepsilon} \left(\frac{\Omega}{\omegap}\right).& 
\end{eqnarray}
From the conservation law (\ref{eq-GP-cons}) we obtain that the wave function $u$ is now normalized to unity:
\begin{equation}
\label{eq-scal-cons}
\| u \|_2^2 = \int_{\R^d} \left|u(\vadim{\vec x})\right|^2\, d\vadim{\vec x} = 1. 
\end{equation}
The total angular momentum (\ref{eq-GP-Lzint}) is now scaled in units of $\hbar$:
\begin{equation}
\label{eq-scal-Lzint}
\vadim{L}_z = \frac{L_z}{N \hbar}= i  \int_{\R^d} u^* \agrad u d\vadim{\vec x}=  \int_{\R^d} Real\left(i u^* \agrad u \right) d\vadim{\vec x}.
\end{equation}

In this non-dimensional setting, the energy (\ref{eq-GP-energ-Agrad}) takes the form:
\begin{equation}
\label{eq-scal-energ-Agrad}
\displaystyle {E}(u) = 
\varepsilon \, \int_{\R^d} \left[ \frac{1}{2} \left|\nabla u + i C_\Omega\, A^t u\right|^2 + C^\text{eff}_\text{trap}\, |u|^2 + \frac{1}{2} C_g |u|^4 \right] \, d\vadim{\vec x},
\end{equation}
and the grand potential (\ref{eq-GP-energ-grandp}) becomes:
\begin{equation}
\label{eq-scal-grandp}
\vadim{\Xi} = E(u) - \vadim{\mu}  \int_{\R^d} |u|^2 \, d\vadim{\vec x}, \quad \vadim{\mu}= \frac{\mu}{\hbar \omega_\perp}.
\end{equation}
The non-dimensional  effective trapping potential corresponding to (\ref{eq-GP-Veff}) is defined as:
\begin{equation}
\label{eq-scal-Ctrap-eff}
C^\text{eff}_\text{trap} =C_\text{trap}-\frac{1}{2} C_\Omega^2 \,{\vadim{r}}^2 =
\frac{1}{\varepsilon^2}\, \left(\vadim{V}(\vadim{x},\vadim{y},\vadim{z}) -\frac{1}{2} \left(\frac{\Omega}{\omegap}\right)^2 \, {\vadim{r}^2}\right) = \frac{1}{\varepsilon^2}\, \vadim{V}^\text{eff}(\vadim{x},\vadim{y},\vadim{z}).
\end{equation}
Finally, the dimensionless form of the stationary GP equation (\ref{eq-GP-stat}) becomes with this scaling:
\begin{equation}
\label{eq-scal-GP-stat}
-\frac{1}{2} \nabla^2 u + C_\text{trap} u + C_g |u|^2 u - i C_\Omega \agrad u = \frac{1}{\varepsilon} \vadim{\mu}\, u,  
\end{equation}

For the trapping potential, we consider in the following a general {\em quadratic+quartic} form that allows to recover the expressions used in most of the theoretical and experimental studies of rotating BEC. Starting from the following physical form of the trapping potential (harmonic potential + detuned laser beam, see \cite{BEC-physV-2004-bretin}):
\begin{equation}
\label{eq-GP-trap-V}
V_\text{trap}(x,y,z)=\frac{m}{2}\left(
\omega_{x}^2 {x}^2 + \omega_y^2 {y}^2 + \omega_z^2 {z}^2 \right)
+ U_2 \left(\frac{r}{w_2}\right)^2  
+ U_4 \left(\frac{r}{w_4}\right)^4,
\end{equation}
we obtain from (\ref{eq-scal-Ctrap}) and (\ref{eq-scal-Ctrap-eff})  the dimensionless effective potential: 
\begin{equation}
\label{eq-scal-trap-V}
\vadim{V}^\text{eff}(\vadim{x},\vadim{y},\vadim{z})=
\frac 12 \left(a_x {\vadim{x}}^2 + a_y {\vadim{y}}^2 + a_z {\vadim{z}}^2
+ a_4 {\vadim{r}}^4\right).
\end{equation}
The non-dimensional coefficients are:
\begin{equation}
\label{eq-scal-trap-ax}
\left\{
\begin{array}{lll}\vspace{0.2cm}
\displaystyle a_x= \left(\frac{\omega_x}{\omegap}\right)^2- 
 \left(\frac{\Omega}{\omegap}\right)^2 +
2\left(\frac{U_2}{m \omegap^2 w_2^2}\right), & &\\ \vspace{0.2cm}
\displaystyle a_y= \left(\frac{\omega_y}{\omegap}\right)^2- 
 \left(\frac{\Omega}{\omegap}\right)^2+ 
2\left(\frac{U_2}{m \omegap^2 w_2^2}\right),&\\ \vspace{0.2cm}
\displaystyle a_z= \left(\frac{\omega_z}{\omegap}\right)^2,&\\ \vspace{0.2cm}
\displaystyle a_4=\frac{2}{\varepsilon}\left(\frac{U_4\, \aho^2}{m \omegap^2 w_4^4}\right) &
\end{array}
\right.
\end{equation}

The classical scaling used in the physical literature is recovered for $\varepsilon=1$. In some mathematical studies \cite{BEC-math-2001-Aftalion-Riviere,BEC-book-2006-Aftalion} it was convenient to define $\varepsilon$ as:
\begin{equation}
\label{eq-scal-AR}
\varepsilon  =\left(\frac{\aho}{8\pi N a_s}\right)^{2/5}. 
\end{equation}
This second scaling, referred as the Aftalion-Rivi{\`e}re (AR) scaling, is particularly appropriate for the Thomas-Fermi (TF) regime characterized by strong interactions (the kinetic energy is negligible compared to the interaction energy). This regime is attained when $N a_s/\aho \gg 1$, which is typically the case in experiments (\eg  \cite{BEC-physV-2000-Madison-a,BEC-physV-2001-Madison,BEC-physV-2002-Rosenbusch-a,BEC-physV-2004-bretin}). In this case, $\varepsilon$ is a small parameter ($\varepsilon \approx 10^{-2}$ in experiments). As a consequence, we notice from (\ref{eq-scal-Cg}) that $C_g   = \frac{1}{2 \varepsilon^2}$ and the GP energy (\ref{eq-scal-energ-Agrad})  becomes:
\begin{equation}
\label{eq-scal-energ-Agrad-TF}
\displaystyle {E}(u) = 
\int_{\R^d} \left[ \frac{\varepsilon}{2} \left|\nabla u + i C_\Omega\, A^t u\right|^2 + \frac{1}{\varepsilon} \vadim{V}^\text{eff}\, |u|^2 + \frac{1}{4 \varepsilon} |u|^4 \right] \, d\vadim{\vec x},
\end{equation}
which is indeed dominated by the trapping and interaction terms. The AR scaling was successfully used in numerical simulation of 2D \cite{BEC-math-2001-Aftalion-Du} or 3D \cite{dan-2003-aft,dan-2004-aft,dan-2005} BEC configurations with vortices.

\newpage

\section{Numerical methods: direct minimisation of the GP energy}\label{section-num}
We present in this section two numerical methods to compute minimizers $u(\vadim{\vec{x}})$ of the non-dimensional GP energy (\ref{eq-scal-energ}) or (\ref{eq-scal-energ-Agrad}), with the constraint (\ref{eq-scal-cons}).  The problem is set on a bounded domain ${\cal D} \in \R^d$, and homogeneous Dirichlet boundary conditions $u=0$ are imposed on $\partial {\cal D}$. The dimensions of ${\cal D}$ will be estimated from the Thomas-Fermi approximation (see Appendix A), in order to ensure that the condensate lies inside ${\cal D}$. The parameters of the minimization problem are the angular velocity $C_\Omega$, the non-linear interaction constant $C_g$ and the trapping potential $C_\text{trap} (\vadim{\vec{x}})$. For the sake of simplicity, the tilde notation for non-dimensional variables will be dropped in the following. 
\subsection{A steepest descent method based on Sobolev gradients}
The first method implemented in our toolbox is the steepest descent method using the Sobolev gradients suggested in \cite{dan-2010-SISC,dan-2010-JCP}.
The algorithm starts from an initial state $u_0(\vec{x})$ and iterates following
\begin{equation}
\label{eq-num-descent}
       u_{n+1} = u_{n} - \alpha_n \, \grdGn,
\end{equation} 
where $\mathcal G_n$ represents the gradient of the energy functional at step $n$ and $\alpha_n$ the descent step. The idea introduced in \cite{dan-2010-SISC} was to define a gradient related to the form (\ref{eq-scal-energ-Agrad}) of the energy. A new Hilbert space, denoted by $H_A({\cal D}, \C)$, was defined and equipped with the inner product:
\begin{equation}
\label{eq-num-pscalHA}
\pscal{u}{v}_{H_A } = \int_{\mathcal D}  \pscal{u}{v}  +  \pscal{\nabla_A u}{\nabla_A v}, 
\end{equation}
where $\nabla_A = \nabla + i C_\Omega A^t$ and $\pscal{u}{v}= u v^*$ denotes the complex inner product. It was proved in \cite{dan-2010-SISC} that the norm  arising from the metric $\| \cdot \|_{H_A }$ is equivalent to the standard Sobolev $H^1$ norm. Hence the completion of $C^1({\mathcal D}, \C)$ with respect to this metric consists of all members of $H^1$.
As a consequence,  the Riez representation theorem in the Hilbert space $H_A = H^1$ allows to define the Sobolev gradient $\nabla_{H_A} E(u)$ as the unique member of $H^1$ such that, $\forall v \in H^1({\cal D}, \C)$:
\begin{equation}
\label{eq-num-Riesz}
E'(u)v= \Real\pscal{\nabla_{L^2} E(u)}{v}_{L^2} = \Real\pscal{\nabla_{H^1} E(u)}{v}_{H^1} = \Real\pscal{\nabla_{H_A} E(u)}{v}_{H_A}. 
\end{equation}  
Since the $L^2$ gradient of the GP energy can be easily derived from (\ref{eq-scal-energ}):
\begin{equation}
\label{eq-num-gradL2}
\nabla_{L^2} E(u)= 2 \varepsilon \, \left( -\frac{1}{2}\nabla^2 u +  C_\text{trap} u +  C_g |u|^2 u -  i C_\Omega \agrad u  \right), 
\end{equation}  
the relationship (\ref{eq-num-Riesz})  allows to compute the $H_A$ gradient. Before using this gradient in the descent method (\ref{eq-num-descent}), it will be projected onto the tangent space of the constraint (\ref{eq-scal-cons}). An explicit projection formula is derived in \cite{dan-2010-SISC}. This technique is an alternative of the usual approach that  re-normalise the solution $u_{n+1}$ after each descent step. 

Compared to the descent method presented in \cite{dan-2010-SISC,dan-2010-JCP}, where a fixed value of the descent step was used, the present method introduces an efficient estimation of the optimal descent step. Since general purpose line-search methods (Brent, Armijo, etc) proved very time consuming for this problem, we finally used the particular line-minimisation analysis specific to the GP energy. The minimiser $\alpha_n$ of the real function:
      \begin{equation} 
      \label{eq-num-Jn}
      J_n(\alpha) = E (u_n - \alpha \, \grdGn),\quad \alpha > 0. 
      \end{equation}
is a root of the third order polynomial:
\begin{equation}
     \label{eq-num-Jnp}
J_n'(\alpha) = c_3 \alpha^3 + c_2 \alpha^2 + c_1 \alpha + c_0,
\end{equation}
with coefficients
\begin{eqnarray}
     \label{eq-num-Jnc1}
c_3 = 2 C_g \int_{\mathcal{D}}  \mod{\grdGn}^4, & \\      \label{eq-num-Jnc2}
c_2 = -6 C_g \int_{\mathcal{D}}   \mod{\grdGn}^2 \, \Real\left(\pscal{u_n}{\grdGn}\right), &\\      \label{eq-num-Jnc3}
c_1 = \int_{\mathcal{D}} \mod{\nabla \grdGn}^2 + 2 C_\text{trap} \mod{\grdGn}^2 + 2 C_g \mod{u_n}^2 \mod{\grdGn}^2 + 4 C_g \Real\left(
\pscal{u_n}{\grdGn}\right)^2 - 2 C_\Omega \Real\left(
i \grdGn^* \agrad \grdGn \right), & \\      \label{eq-num-Jnc4}
c_0 = -\int_{\mathcal{D}} \Real\left(\pscal{\nabla u_n}{\nabla \grdGn}\right) + 2 \Real\left(
\pscal{u_n}{\grdGn}\right) \left[C_\text{trap} + C_g \mod{u_n}^2\right] - 2 C_\Omega \Real\left(
i \grdGn^* \agrad u_n \right). &
\end{eqnarray}
In FreeFem++, we can use the function \textbf{polycomplexsolve} (from \textbf{GSL} library)  \cite{gsl} to  calculate the three  roots of  the polynomial $J_n'(\alpha)$ and then select the root realizing the minimum of the  energy $J_n(\alpha)$.

The algorithm for the descent method can be easily identified in the FreeFem++ scripts, since appropriate macros were defined for the mathematical operators (inner product, norms, etc). All variables are discretised using $P^1$ finite elements; the non-linear term is represented with $P^4$ finite elements in 2D and $P^2$ in 3D. The following steps were programmed, with a syntax very close to mathematical relationships:
\begin{enumerate}
\item Suppose that  the solution $u_n$ at iteration $n$ was built.
We compute $G = \nabla_{H_A}E(u_n)/(2\varepsilon)$, solution of the variational problem corresponding to (\ref{eq-num-Riesz}) and (\ref{eq-num-gradL2}):

$\forall v \in {H_0^1(\mathcal D, \mathbb C)},$

\begin{equation}
     \label{eq-num-algoS1}
\int_{\mathcal D} \left( 1 + C_\Omega^2 (x^2 + y^2) \right) G v + \nabla G \nabla v - 2i C_\Omega A^t \nabla G v = \int_{\mathcal D} \frac{1}{2} \nabla u_n \nabla v + \left[ C_\text{trap} u_n + C_g \lvert u_n \rvert ^2 u_n - iC_\Omega A^t \nabla u_n \right] v.
\end{equation}
\item We compute the projection of G over the tangent space of the unitary norm constraint (see \cite{dan-2010-SISC}):
\begin{equation}
     \label{eq-num-algoS3}
 P_{u_n,H_A}G = G - \frac{ \Real \left( \langle u_n , G \rangle_{L^2} \right) }{ \Real \left( \langle u_n , v_{H_A} \rangle_{L^2} \right) } v_{H_A},
\end{equation}
where $v_{H_A}$ is solution of the variational problem:
\begin{equation}
     \label{eq-num-algoS3b}
 \langle v_{H_A} , v \rangle_{H_A} = \langle u , v \rangle_{L^2}, \ \forall v \in {H_A}.
\end{equation}

\item We compute the optimal descent step:
\begin{equation}
     \label{eq-num-algoS4}
 \chi_n = \min_{\alpha >0} \, E (u_n - \chi \, P_{u_n,H_A}G),
\end{equation}
by finding the roots of the third order polynomial (\ref{eq-num-Jnp}) with coefficients (\ref{eq-num-Jnc1})-(\ref{eq-num-Jnc4}) and choosing the one realizing the minimum of the line energy. Note that the factor $(2\varepsilon)$ appearing in the expression of the gradient (\ref{eq-num-gradL2}) was included in the expression of the optimal descent step $\chi = (2\varepsilon)\alpha$.

\item We build the solution $u_{n+1}$ at iteration $n+1$:
\begin{equation*}
    \label{eq-num-algoS5}
 u_{n+1} = u_n -  \chi_n \  P_{u_n,H_A}G.
\end{equation*}

\item Finally, we compute the relative error  $\delta E_{n+1} = \frac{E(u_{n+1})-E(u_n)}{E(u_{n+1})}$ and call the mesh adaptivity algorithm suggested in \cite{dan-2010-JCP} (see below). Convergence to the stationary state is achieved if $\delta E_{n+1} < \varepsilon_c=10^{-9}$.
\end{enumerate}

\subsection{Mesh adaptation}\label{sec-mesh-adapt}
FreeFem++ includes a powerful  mesh adaptivity tool (function \textbf{adaptmesh}) using metric control algorithms suggested in  \cite{hecht-1996-missi,hecht-2000-ijnmf,hecht-1997-aiaa,george-1998}. The main idea is to define a metric based on the Hessian and use a Delaunay procedure to build a new mesh such that all the edges are close to the unit length with respect to this new metric. In the steepest descent algorithm, we call this function after building $u_{n+1}$ in the step 4 of the previous algorithm. Since our convergence criterion is based on the relative change of energy of the solution  ($\delta E_{n+1}$) we use the same indicator to trigger the mesh adaptive procedure following the next algorithm (see also \cite{dan-2010-JCP}):
\begin{enumerate}
  \item choose a sequence of decreasing values $\varepsilon^i \ge \varepsilon_c$, that represent threshold values for the mesh adaptivity;
  \item set $i=1$;
  \item if $\delta E_{n+1}$ is decreasing and $\varepsilon^{i+1} < \delta E_{n+1} < \varepsilon^{i}$ and $\delta E_n > \varepsilon_c$, call the mesh adaptivity procedure; the solution $u$ is interpolated on the new mesh and normalized to satisfy the unitary norm constraint;
  \item if $\delta E_{n+1}$ is increasing, \ie large variations of the energy appear (\eg if new vortices enter the domain), reconsider the previous bounds by setting $i \rightarrow {i-1}$;
  \item if step 3 was performed $N_{ad} \geq 1$ times,  increase $i$ to $i+1$. Limiting the number of mesh refinements for the same threshold, is necessary since, at step 2, the interpolation on the new refined mesh and the normalization of the solution could lead to an increase of the value of $\delta E_{n+1}$.
\end{enumerate}

\begin{figure}[h]
\begin{center}
	\includegraphics[width=0.65\textwidth]{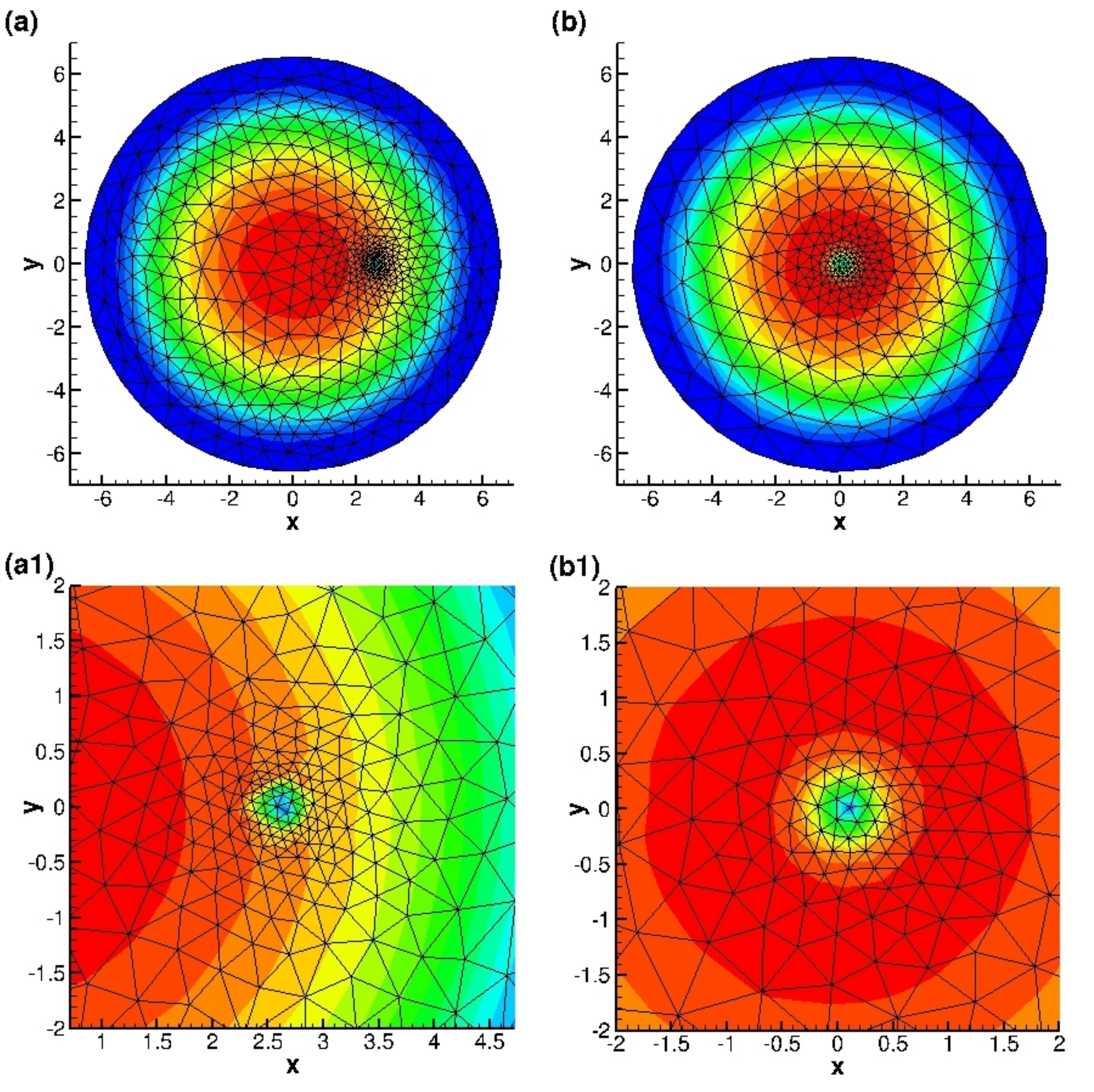}
	\end{center}
\caption{Illustration of the mesh adaptivity in 2D. Test case with an initial state containing an off-centred vortex (a) and a final (stationary) state with a central vortex (b). The mesh refinement follows the evolution of the vortex position (corresponding zoom in figures a1 and b1).}
\label{fig-2D-adapt}
\end{figure}

Figure \ref{fig-2D-adapt} illustrates the 2D mesh adaptivity procedure. It represents a test case where the initial field has an off-centred vortex and  the final (converged) solution is expected to present a centred vortex (the details of the parameters for this case are given in section \ref{sec-example}). We plot in figure \ref{fig-2D-adapt}(a) the initial state, built with the Thomas-Fermi approximation. In figure \ref{fig-2D-adapt}(b) we plot the final solution converged with the Sobolev gradient method. The corresponding zoom in the vortex area are displayed in figures \ref{fig-2D-adapt} (a1) and (b1). Note that the mesh adaptivity procedure generated a denser mesh near the position of the vortex. Similarly, the number of triangles was decreased near the border, where the solution is smoother. However, the de-refinement of the mesh must be used with caution for high rotation rates, when new vortices can nucleate in the condensate near the boundary. It was shown in \cite{dan-2010-JCP} that this procedure decreases the CPU time and the number of iterations when compared with computation on  fixed refined meshes. 

\begin{figure}[h]
\begin{center}
\includegraphics[width=0.45\textwidth]{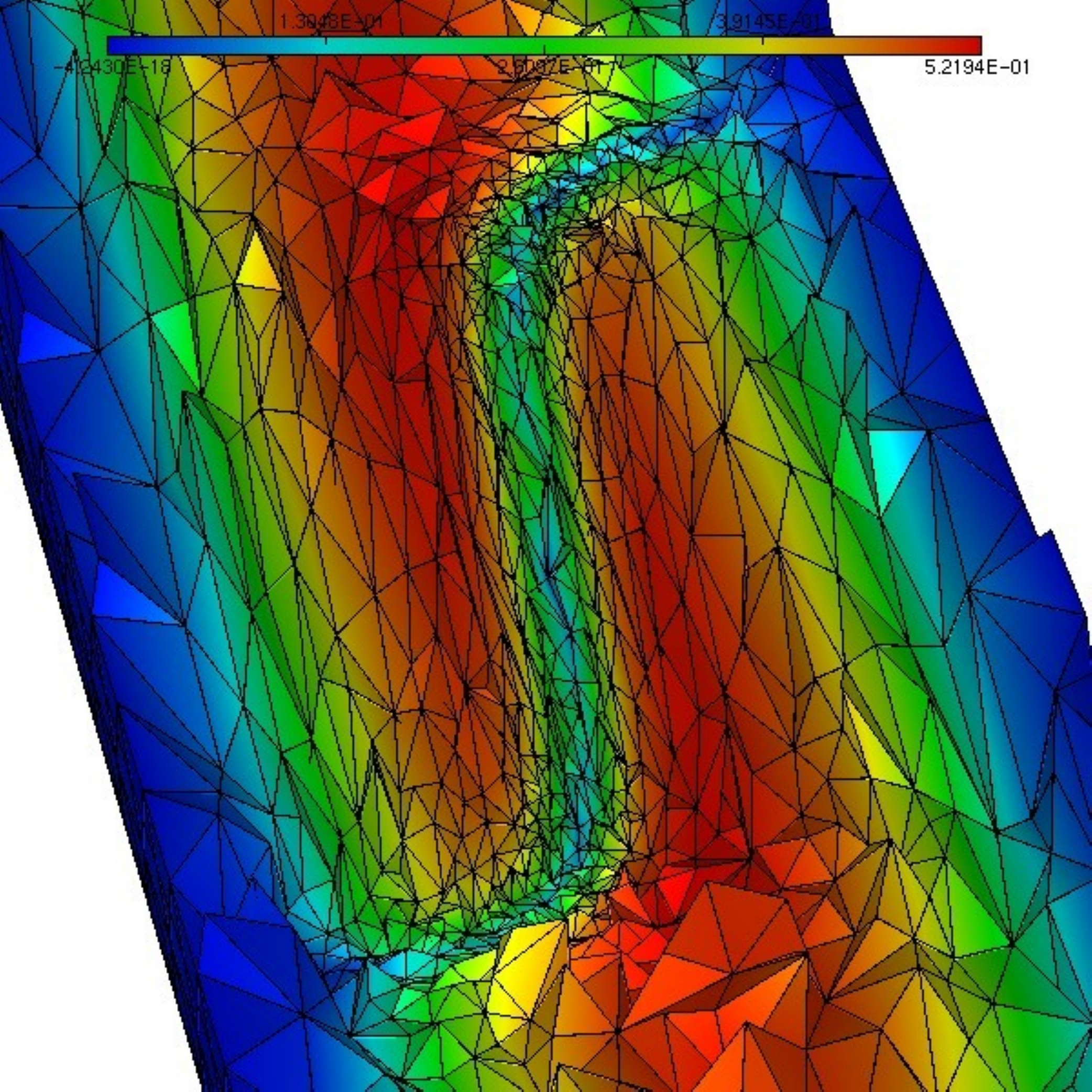}
\end{center}
\caption{Illustration of the mesh adaptivity in 3D. Test case computing the equilibrium configuration with a single S-shape vortex line.}
\label{fig-3D-adapt}
\end{figure}

For  3D computations, FreeFem++ uses the function \textbf{mshmet} \cite{mshmet} to compute the metrics and the function \textbf{mmg3d} \cite{mmg3d} to build the new mesh corresponding to this metric. In figure \ref{fig-3D-adapt}, we plot a 3D mesh adapted to the solution presenting a vortex line with a "S" shape. We carried out the visualisation with \textbf{medit}, a mesh visualisation software \cite{medit} interfaced with FreeFem++.  Note that the mesh adaptation follows precisely the vortex line by adding triangles for a better accuracy. Outside the vortex area, the mesh adaptation allowed us to have less triangles, with a bigger size. 

\pagebreak\clearpage

\subsection{Minimisation algorithm using the optimisation library Ipopt}
The optimisation library Ipopt is based on an interior point minimisation method \cite{Ipopt-int-point-2002}, a barrier functions tool  \cite{Ipopt-barr-strat-2008} and a filter line search \cite{Ipopt-line-search-2005}. This powerful state-of-the-art optimisation library is interfaced with FreeFem++ \cite{Auliac-thesis-2014} and offers the possibility to solve constrained optimisation problems of the general form:
\begin{equation}
\label{eq-ipopt-Prob}
 \text{find } {\bf x_0} = \underset{{\bf x} \in \mathbb{R}^n}{\text{argmin}} (f({\bf x})),
\end{equation}
\begin{equation}
\label{eq-ipopt-Probc}
\text{ such that } \left\{\begin{array}{ll} 
\forall i \leq n, \ x^{lb}_i \leq x_i \leq x^{ub}_i \text{ (simple bounds)},\\ 
\forall i \leq m, \ c^{lb}_i \leq c_i({\bf x}) \leq c^{ub}_i \text{ (constraint functions)},
\end{array}\right.
\end{equation}
where $lb$ stands for {\em lower bound} and $ub$ for {\em upper bound}. If for some $i \leq m, \ c^{lb}_i = c^{ub}_i$ we obtain an equality constraint.

For the minimisation of the Gross-Pitaevskii energy, the use of Ipopt is quite simple: the conservation constraint (\ref{eq-scal-cons}) is an 
equality constraint and, consequently, we take $m = 1$ and $c^{lb} = c^{ub} =1$ in the previous general form. Ipopt will then solve the Euler-Lagrange equation associated to the problem (\ref{eq-ipopt-Prob})-(\ref{eq-ipopt-Probc}):
\begin{equation}
\left\{\begin{array}{ll} 
\nabla f({\bf x}) + \lambda \, \nabla c({\bf x}) = 0,\\ 
c({\bf x}) = 0,
\end{array}\right.
\end{equation}
where $\lambda \in \mathbb R$ is a Lagrange multiplier and $c({\bf x})$ the constraint. Note that, in our case, $\lambda$ corresponds to the chemical potential.
Let us define 
\begin{equation}
L({\bf x}, \lambda) := f({\bf x}) + \lambda \, c({\bf x}).
\end{equation}
Ipopt first finds a descent direction $({\bf dx} , \, d\lambda)$ by using the Newton method. Indeed, at each iteration $n$ it solves the system: 
\begin{equation}
\begin{pmatrix}
\nabla^2 L({\bf x_n}, \lambda_n) \quad \nabla c({\bf x_n})\\
\nabla c({\bf x_n}) \quad 0
\end{pmatrix}
\begin{pmatrix}
{\bf dx}\\
d\lambda
\end{pmatrix}
= -
\begin{pmatrix}
\nabla L({\bf x_n}, \lambda_n)\\
c({\bf x_n})
\end{pmatrix}
\end{equation}
Then it advances at the next step:
\begin{equation*}
\begin{pmatrix}
{\bf x_{n+1}}\\
\lambda_{n+1}
\end{pmatrix}
= \begin{pmatrix}
{\bf x_n}\\
\lambda_n
\end{pmatrix}
+ \alpha_n \,
\begin{pmatrix}
{\bf dx}\\
d\lambda
\end{pmatrix}
\end{equation*}
where $\alpha_n \in \, (0,1]$ is a descent step computed using the filter line-search method suggested in \cite{Ipopt-line-search-2005}.
The algorithm will stop when either the error $(\varepsilon_n = \max \, (\|\nabla f({\bf x_n}) + \lambda_n \, \nabla c({\bf x_n})\|_\infty, \, \|c({\bf x_n})\|_\infty )\,)$ or the number of iterations reaches a value defined by the user.

As Ipopt seeks for solutions in $\mathbb R^n$, we have to separate in the  Gross-Pitaevskii energy functional the real and imaginary part.
The problem to solve becomes:
\begin{equation*}
\text{find } [u_r,u_i] \in (H_0^1(\mathcal D,\mathbb R))^2 \text{ wich minimizes}  
\end{equation*}
\begin{equation}
 E(u_r,u_i) =  \int_{\mathcal{D}} \left[ \frac{1}{2}\, \lvert \nabla u_r \rvert^2 + \frac{1}{2}\, \lvert \nabla u_i \rvert^2 + C_{trap}\, (u_r^2  +u_i^2) + \frac{1}{2} \,C_g\, ( u_r^2  +u_i^2 )^2 \right]\\
 -\,C_{\Omega}\, L_z(u_r,u_i),
\end{equation}
with
\begin{equation}
L_z(u_r,u_i) = \int_{\mathcal{D}} \left[  y \,\left( \frac{\partial u_r}{\partial x}\,u_i - \frac{\partial u_i}{\partial x}\,u_r\right) - x \,\left(   \frac{\partial u_r}{\partial y}\,u_i - \frac{\partial u_i}{\partial y}\,u_r \right) \right].
\end{equation}
Then we can calculate the Frechet derivative of $E$ as:
\begin{align}\nonumber
 E'(u_r,u_i)\cdot[v_r,v_i] &=  \int_{\mathcal{D}} \left[ \nabla u_r \cdot  \nabla v_r + \nabla u_i \cdot  \nabla v_i
      + 2 \,C_{trap} \, (u_r\,v_r + u_i\,v_i) \right]\\ \nonumber
      &+ 2 \,C_g \,\int_{\mathcal{D}}( u_r^2 + u_i^2 ) \, (u_r\,v_r + u_i\,v_i)\\
      &- C_{\Omega} L_z'(u_r,u_i)\cdot[v_r,v_i],
\end{align}
with
\begin{align}  \nonumber
 - L_z'(u_r,u_i)\cdot[v_r,v_i] = & \int_{\mathcal{D}}  y \,\left[ - \frac{\partial u_r}{\partial x}\,v_i +\frac{\partial u_i}{\partial x}\,v_r - \frac{\partial v_r}{\partial x}\,u_i +\frac{\partial v_i}{\partial x}\,u_r\right]\\
 + & \int_{\mathcal{D}} x \,\left[ \frac{\partial u_r}{\partial y}\,v_i - \frac{\partial u_i}{\partial y}\,v_r + \frac{\partial v_r}{\partial y}\,u_i - \frac{\partial v_i}{\partial y}\,u_r  \right].
\end{align}
Finally, the calculus of the the second order Frechet derivative of of E leads to:
\begin{align}\nonumber
  E''(u_r,u_i)\cdot([v_r,v_i],[w_r,w_i]) &=  \int_{\mathcal{D}} \left[ \nabla v_r \cdot  \nabla w_r + \nabla v_i \cdot  \nabla w_i
      + 2 \,C_{trap} \, (v_r\,w_r + v_i\,w_i) \right]\\ \nonumber
      &+ 2 \, C_g \,\int_{\mathcal{D}}\left[( u_r^2 + u_i^2 ) \, (v_r\,w_r + v_i\,w_i)\right]\\ \nonumber
      &+ 4 \, C_g \,\int_{\mathcal{D}}\left[(u_r\,v_r + u_i\,v_i) \, (u_r\,w_r + u_i\,w_i)\right]\\
      &- C_{\Omega} L_z''(u_r,u_i)\cdot([v_r,v_i],[w_r,w_i]),
\end{align}
with
\begin{align} \nonumber 
   - L_z''(u_r,u_i)\cdot([v_r,v_i],[w_r,w_i]) = & \int_{\mathcal{D}}  y \,\left[- \frac{\partial w_r}{\partial x}\,v_i +\frac{\partial w_i}{\partial x}\,v_r - \frac{\partial v_r}{\partial x}\,w_i +\frac{\partial v_i}{\partial x}\,w_r\right]\\
 + & \int_{\mathcal{D}} x \,\left[ \frac{\partial w_r}{\partial y}\,v_i - \frac{\partial w_i}{\partial y}\,v_r + \frac{\partial v_r}{\partial y}\,w_i - \frac{\partial v_i}{\partial y}\,w_r  \right].
\end{align}
The expression of the constraint functional is:
\begin{equation}  
 c(u_r,u_i) =  \int_{\mathcal{D}} ( u_r ^2 + v_r ^2 ) - 1,
\end{equation}
and its gradient:
\begin{equation}  
 \nabla c(u_r,u_i)\cdot[v_r,v_i] = 2 \, \int_{\mathcal{D}} (u_r \, v_r + u_i \, v_i).
\end{equation}

With Ipopt linked as an external library to FreeFem++, we can not directly use mesh adaptivity in its internal algorithm. In exchange, we can couple the computation of the minimizer with the mesh adaptivity procedure.  The following algorithm was implemented in the programs.
Set $n_\text{adapt}$, the total number of mesh refinements to be done and $\varepsilon_0$ and $\varepsilon_\text{last}$, the first and the last mesh adaptivity prescribed errors (parameters of the FreeFem++ function \textbf{adaptmesh}).
\begin{enumerate}
\item At step $k \in [0, n_\text{adapt}-1]$, run Ipopt to find a solution $[u_r^k,u_i^k]$.
\item Build a new mesh adapted to $[u_r^k,u_i^k]$ with a prescribed mesh adaptivity error
 \begin{equation}\label{eq-algo-ipopt}
\varepsilon_k = \varepsilon_0 \, \left( \frac{\varepsilon_\text{last}}{\varepsilon_0} \right)^{k/(n_\text{adapt}-1)}.
\end{equation}
\item Go to step $k+1$.
\end{enumerate}
Typical values used  for 2D computations are $n_\text{adapt} = 4$, $\varepsilon_0= 0.1$ and $\varepsilon_\text{last}= 0.005$. For 3D cases, as the computation is more difficult, it's more convenient to use a higher number of mesh adaptations and a lower ratio $\varepsilon_\text{last}/\varepsilon_0$. Typical values are $n_\text{adapt} = 6$, $\varepsilon_0= 0.01$ and $\varepsilon_\text{last}= 0.005$.

\section{Building the initial approximation}

In computing stationary states for rotating BEC, the initial approximation used to start the iterative methods is of crucial importance. It can not only affect the convergence speed, but also the topology of the stationary solution, especially when local minima (meta-stable) solutions are sought. 
We present in this section three methods to build initial states for the computation of stationary solutions: the Thomas-Fermi approximation, a rapid calculation of the ground state with Ipopt for simplified configurations (axisymmetric or non-rotating) and, finally, an ansatz for a manufactured initial state with vortices.

\subsection{Analytical solution based on the Thomas-Fermi approximation}\label{sec-TF}
The Thomas-Fermi regime is characterized by strong interactions (the kinetic energy is negligible compared to the interaction energy). This regime is attained when $N a_s/\aho \gg 1$. 
If the healing length
$\xi=\left(\displaystyle 8\pi a_s \rho\right)^{-1/2}$, with  $\rho$ the atomic density,
is defined as the length for which the kinetic and interaction energies are comparable, in the Thomas-Fermi regime the characteristic length scales are larger that the healing length. We give below some typical values from experiments of BEC with vortices \cite{BEC-physV-2000-Madison-a,BEC-physV-2000-Madison-b,BEC-physV-2003-bretin-prl}:
\begin{equation*}
\begin{array}{ccccccccc}
a_s & \ll & 1/\rho^{1/3} & < & \xi & \ll & \aho & \ll & R \\
5\, [nm] & \ll & 0.2 \, [\mu m] & < & 0.3 \, [\mu m]  & \ll & 1 \, [\mu m] & \ll & 5 \, [\mu m],
\end{array}
\end{equation*}
where $1/\rho^{1/3}$ approximates the distance between atoms and $R$ is the radius of the condensate.

The general form of the Thomas-Fermi approximation of the atomic density ($\rho = |u|^2$) is obtained by neglecting the first term in the energy (\ref{eq-scal-energ-Agrad}). The Euler-Lagrange equation of the corresponding grand potential (\ref{eq-scal-grandp}) gives:
\begin{equation}
\label{eq-scal-TF}
\rtf = \left(\frac{{\vadim{\mu}}/\varepsilon - C^\text{eff}_\text{trap}}{C_g}\right)_+  =  \frac{1}{\varepsilon^2 C_g} \left( \varepsilon \frac{\mu}{\hbar \omega_\perp} -  \vadim{V}^\text{eff}\right)_+.
\end{equation} 
We notice that this form is equivalent to the usual Thomas-Fermi approximation for non-rotating condensates, but with a trapping potential  (\ref{eq-scal-Ctrap-eff}) corrected by the centrifugal term (see also \cite{BEC-phys-1999-Stringari-TF}). Following (\ref{eq-GP-Veff}), for a harmonic trapping potential the radial trapping frequency $\omega_\perp$ is thus replaced by $(\omega_\perp^2 - \Omega^2)^{1/2}$.

It is also interesting to note from (\ref{eq-scal-u}) and (\ref{eq-scal-TF}) that the atomic density in numerical simulations using the AR scaling with typical value $\varepsilon=10^{-2}$ is amplified by a factor of $10^4$, when compared to the classical scaling ($\varepsilon=1$). This remark is important for setting the numerical value which will serve to identify a quantized vortex: since theoretically $\rho=0$ in the vortex centre, the low value $\rho_\text{min}$ of the iso-contour level used to represent vortices will depend on the scaling.

We use in the following the Thomas-Fermi approximation to estimate the size of the computational domain and also to set the initial guess for the minimisation algorithms. We derive in Appendix A closed formulae for the Thomas Fermi approximation corresponding to different types of potentials: harmonic, {\em quartic+quadratic}, {\em quartic-quadratic}.

\subsection{Numerical approximation with Ipopt for axisymmetric or non-rotating cases}

The main drawback of the Thomas-Fermi approximation, which is generally a truncated parabola, is the discontinuity of its first derivative on the border of the condensate where $\rtf =0$. This could trigger oscillations of the solution, when high-order (spectral) methods are used for the space discretisation. A smoother initial field can be obtained by directly computing with Ipopt a minimizer of the GP energy. When simplified forms of the energy (\eg axisymmetric) are used, this preliminary computation is very cheap in terms of computational time. 

We present below the approach of computing axisymmetric initial fields with Ipopt, corresponding to the ground state without vortices or with a central vortex of given winding number. We consider the cylindrical coordinates $(r, \, \theta, \, z)$ and  assume that the  solution is axisymmetric ($\frac{\partial u}{\partial \theta} = 0$) and symmetric in the $z$-direction ($u(z) = u(-z)$). This  is also the case of the Thomas-Fermi approximation if $a_x=a_y$ in the trapping potential (\ref{eq-scal-trap-V}).

\begin{figure}
\begin{center}
\includegraphics[width=\textwidth]{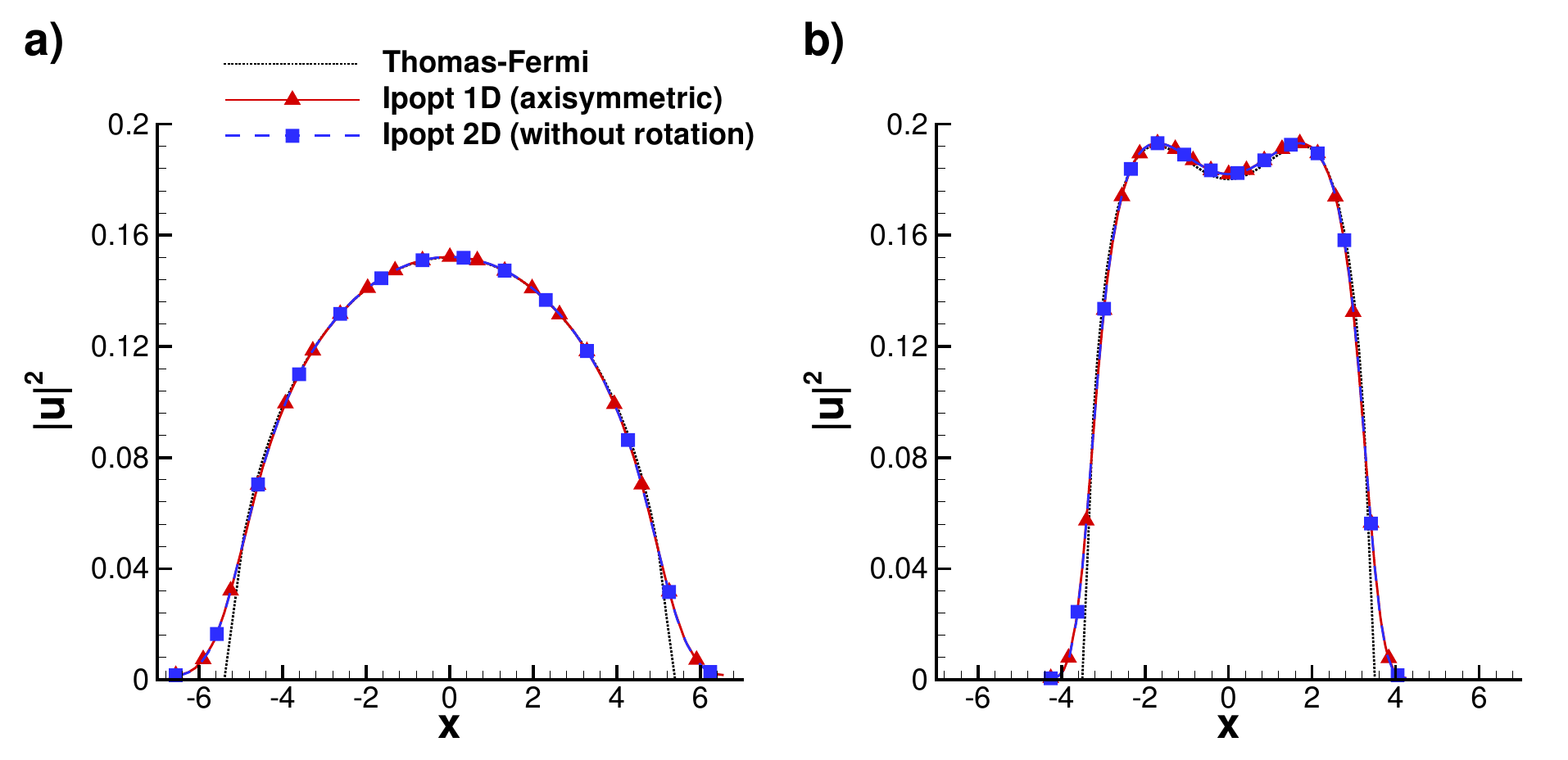}
\end{center}
\caption{Initialisation of a 2D calculation. Density profiles corresponding to the Thomas-Fermi approximation (solid line), the axisymmetric solution computed with Ipopt (\Red{$\blacktriangle$}) and the full 2D solution computed with Ipopt without rotation (\Blue{$\blacksquare$}). Harmonic potential (a) and quartic potential (b).}
\label{figInitComp}
\end{figure}

Since $\frac{\partial u}{\partial \theta} = x\frac{\partial u}{\partial y} - y\frac{\partial u}{\partial x}$ we can infer that $L_z = 0$.
Then the energy becomes 
 \begin{align}\nonumber
 E(u) &= \int_\mathcal D \left[ \frac{1}{2} \left( \mod{\frac{\partial u}{\partial r}}^2 + \mod{\frac1r \frac{\partial u}{\partial \theta}}^2 + \mod{\frac{\partial u}{\partial z}}^2 \right)
      + C^\text{eff}_\text{trap} \lvert u \rvert^2 
      + \frac{1}{2} C_g \lvert u \rvert^4 \right] r dr \, d\theta \, dz.\\     
      &= 4 \pi \int_0^{R_{max}}\int_0^{z_{max}}  \left[ \frac{1}{2} \left(\mod{\frac{\partial u}{\partial r}}^2 + \mod{\frac{\partial u}{\partial z}}^2 \right)
      + C^\text{eff}_\text{trap} \lvert u \rvert^2 
      + \frac{1}{2} C_g \lvert u \rvert^4 \right] r dr \, dz.  
\end{align}
The 3D problem is now reduced to a 2D problem.
In order to solve this 2D problem with Ipopt, we need the Fr\'{e}chet derivative of $E$ and its Hessian: 
\begin{equation}
 E'(u).v = 4 \pi \int_0^{R_{max}}\int_0^{z_{max}}  \left[ \frac{\partial u}{\partial r}   \frac{\partial v}{\partial r} + \frac{\partial u}{\partial z}  \frac{\partial v}{\partial z}
         + 2 C^\text{eff}_\text{trap}   u v 
         + 2 C_g \lvert u \rvert^2 u v \right] r dr \, dz,  
\end{equation}
\begin{equation}
 E''(u)[v,w] =  4 \pi \int_0^{R_{max}}\int_0^{z_{max}} \left[ \frac{\partial v}{\partial r}   \frac{\partial w}{\partial r} + \frac{\partial v}{\partial z}   \frac{\partial w}{\partial z}
      + 2 C^\text{eff}_\text{trap}  v w 
      + 2 C_g \left(v w \lvert u \rvert^2 + 2 u v \Real(uw) \right) \right] r dr dz.
\end{equation}
In the case of a 2D simulation, the axisimmetry reduces the 2D problem to a 1D problem. In this case, the integration and the derivative with respect to $z$ must be omitted in previous formulations.
Figure \ref{figInitComp} shows a comparison between the Thomas-Fermi approximation and the axisymmetric solution computed with Ipopt for two trapping potentials (\ref{eq-scal-trap-V}): harmonic potential, with  $a_x=a_y=1, a_4=0$ and quartic potential , with  $a_x=a_y=1, a_4=0.5$. A third solution, obtained by using the full 2D formulation of the GP energy without rotation in Ipopt, is also plotted for reference. We notice the regularity of the axisymmetric solution in the vicinity of $\rtf=0$ and the good approximation it offers, when compared with the full 2D computation. For anisotropic potentials, we can still use the full (2D or 3D) formulation of the GP energy without rotation to compute with Ipopt an initial condition for the computations with rotation.

\subsection{Manufactured initial state with vortices}\label{sec-manuf}
Sometimes it is necessary to manufacture initial states by artificially including vortices. This could be useful when local minima, corresponding to meta-stable solutions, are sought. If $u(x,y,z)$ is the ground state without rotation (set by the TF approximation or computed with Ipopt), we can add vortices by multiplying  $u$ in each plane ($x,y)$ by the following ansatz used in \cite{dan-2004-aft,dan-2005} for 3D simulations:
\begin{equation}\label{eq-anzatz}
 u_v(x,y) = \sqrt{\frac{1}{2}\left[1+tanh\left(\frac{4}{\varepsilon_v}(r_v-\varepsilon_v)\right)\right]}\cdot e^{i\theta_v}, 
\end{equation}
where $r_v =\sqrt{(x-x_c)^2 + (y-y_c)^2}$ and $\theta_v = atan \left(\frac{y-y_c}{x-x_c}\right)$ are the polar coordinates taken from the imposed centre $(x_c, y_c)$ of the vortex and $\varepsilon_v$ the vortex  radius.
In order to obtain a particular 3D shape of the vortex (U-shaped or S-shaped vortex, see figure \ref{fig-3D-example}), we can prescribe the position of the vortex centre in each transverse plane $(x,y)$. For example, a S-vortex lying in the major $(x,z)$ plane will have $y_c=0$ and 
\begin{equation}\label{eq-anzatz-shape}
\left\{
\begin{array}{lll}\vspace{0.2cm}
\displaystyle x_c(z) = -1+\frac{\tanh\left[\alpha_v \left(1+\frac{z}{\beta_v}\right)\right]}{\tanh(\alpha_v)}, \, \text{  if } z < 0,\\
\displaystyle x_c(z) = 1+\frac{\tanh\left[\alpha_v \left(-1+\frac{z}{\beta_v}\right)\right]}{\tanh(\alpha_v)}, \, \text{  if } z \geq 0,\\
\end{array}
\right.
\end{equation}
where $\alpha_v$ and $\beta_v$ respectively control the curvature and the length of the vortex.

\section{Description of the programs}

The methods described previously were implemented in a 2D and a 3D toolbox based on the FreeFem++ software \cite{freefem}. Using two input files, the toolbox offers to the user the choice between two scalings (classical or Aftalion-Riviere), three ways of computing the initial approximation (Thomas-Fermi, axisymmetric or non-rotating)  and two methods to compute the ground state (Sobolev gradient or Ipopt). The main difference between the 2D and 3D codes is in the post-treatment part: we can automatically count the number of vortices in 2D while it is more difficult in 3D. Also, the setting of input parameters is different: there are several additional parameters in 3D to control the shape of the vortex ansatz (I-shaped, S-shaped or U-shaped) and the shape of the initial mesh (cylindrical or ellipsoidal). Moreover, the user can choose to plot the evolution of the energy with \textbf{Gnuplot} \cite{gnuplot} during the computation and the evolution of the solution with either the FreeFem++ plotting tool or using \textbf{Medit} \cite{medit}. In this section we first describe the architecture of the programs and the organisation of the files. Then we focus on the list of input parameters and the structure of output files.
\subsection{Program architecture}
\begin{figure}[ht]
\begin{center}
\includegraphics[width=\textwidth]{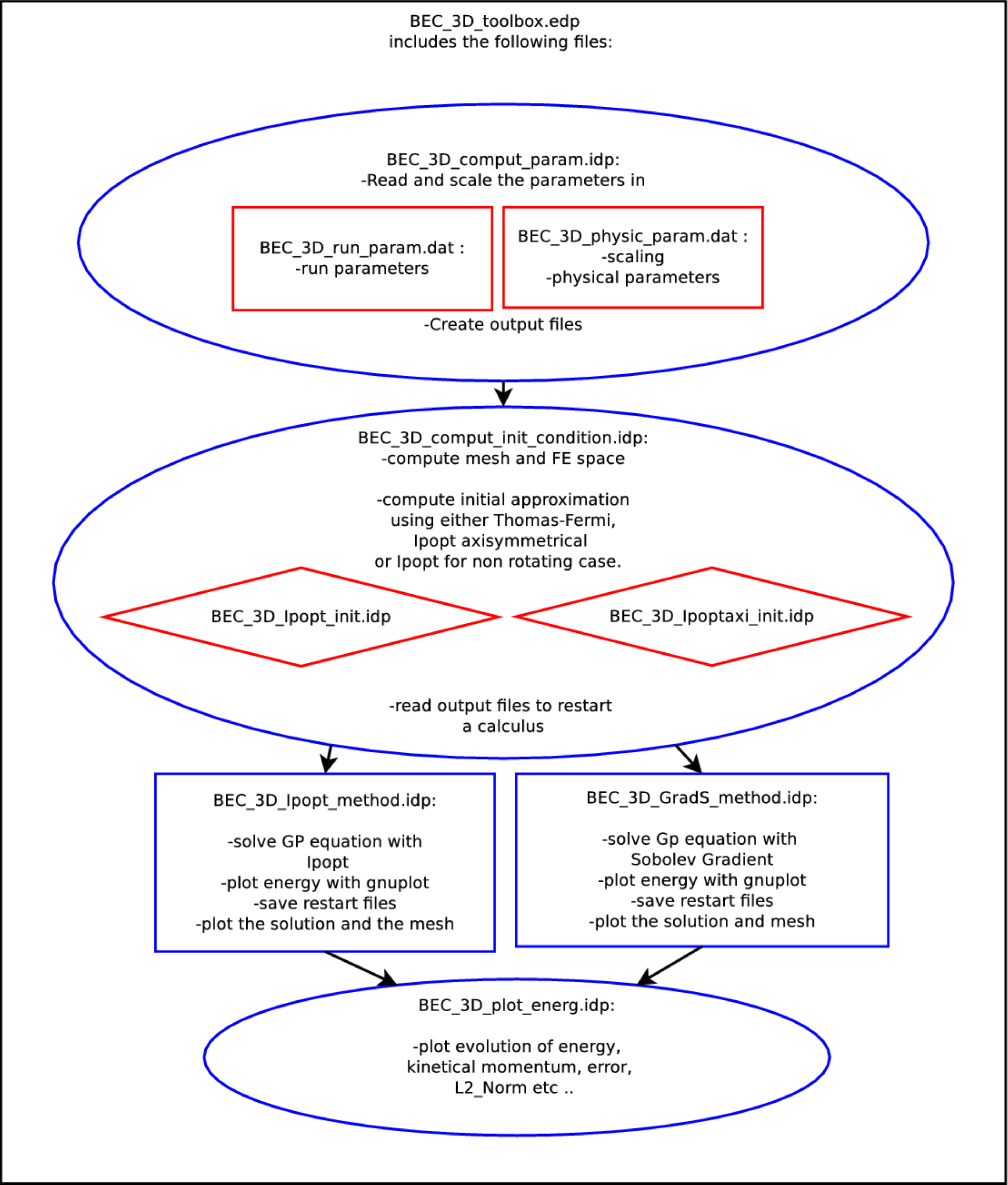}
\end{center}

\caption{Program architecture of the 3D toolbox.}
\label{figDiag}
\end{figure}

Figure \ref{figDiag} gives a schematic overview of the content of the 3D toolbox. The 2D toolbox has similar architecture.
\clearpage
All files are provided in a directory called \textit{BEC\_$\mathcal X$D\_ToolBox\_FreeFem} where $\mathcal X$ is the dimension 2 or 3. This directory contains: 
\begin{enumerate}
 \item The \textit{BEC\_$\mathcal X$D\_ToolBox.edp} file containing the main script.
 \item The \textit{Input} directory where two files allow the user to choose parameters:
 \begin{itemize}
  \item  \textit{BEC\_$\mathcal X$D\_physic\_param.dat} contains the parameters describing the physical case.
  \item  \textit{BEC\_$\mathcal X$D\_run\_param.dat} contain choices for the run.
 \end{itemize}
 \item The \textit{Include} directory which contains 9 files:
 \begin{itemize}
   \item  \textit{BEC\_$\mathcal X$D\_Macros.idp} contains all the usefull macros and functions.
  \item  \textit{BEC\_$\mathcal X$D\_comput\_param.idp} reads the parameters files and builds constants.
  \item  \textit{BEC\_$\mathcal X$D\_comput\_init\_condition.idp} compute an initial approximation using  either Thomas-Fermi or Ipopt.
  \item  \textit{BEC\_$\mathcal X$D\_Ipoptaxi\_init.idp} contains the script to use with Ipopt axisymmetric in dimension ($\mathcal X-1$) for the initial condition.
  \item  \textit{BEC\_$\mathcal X$D\_Ipopt\_init.idp} contains the script to use Ipopt in dimension $\mathcal X$ without rotation to build the initial condition.
  \item  \textit{BEC\_$\mathcal X$D\_GradS\_method.idp} solves the main problem with Sobolev gradient method.
  \item  \textit{BEC\_$\mathcal X$D\_Ipopt\_method.idp} solves the main problem with Ipopt method.
  \item  \textit{BEC\_$\mathcal X$D\_plot\_energ.idp} builds a gnuplot script and runs gnuplot in order to plot the energy and other relevant quantities.
  \item  \textit{BEC\_2D\_results.idp} finds the number of vortices and gives their position in 2D.
 \end{itemize}
 \item The \textit{Examples} directory. In 2D, this directory contains 8 examples of input files allowing the user choose between two cases of scaling, potential or method.
 To do so, in a terminal the user can write, for example, the command line:\\
  FreeFem++ \textit{BEC\_$\mathcal X$D\_ToolBox.edp}\\
 -run \textit{Examples/GradS\_Harm\_run\_param.dat}\\ 
 -param \textit{Examples/AR\_Harm\_physic\_param.dat}.\\
  This will run the program with an harmonic potential, the Aftalion-Riviere scaling and Sobolev gradient method.
 In 3D, this directory contains 6 files to run examples to compute a S-shaped or a U-shaped vortex, using either the Sobolev Gradient method or Ipopt. A more precise description of these examples is provided in section \ref{sec-example}.
 \item \textit{GLUT} directory contains a C++ script that must be compiled to create a user interface with GLUT.
 \item A makefile to compile the source code for the interface and a README file.
\end{enumerate}
\subsection{Input parameters}
We focus now on the description of the input parameters. These are distributed in two files. In both files, comments are preceded by the usual // symbol and key words by the @ symbol. If the user wants to set a parameter, he has to enter its value after the corresponding key word. If a key word is not written in a file, a default value is given to the corresponding parameter. Some parameters must be specified by the user, otherwise the computation stops (see bellow). By default, the user has to use the two files provided in the \textit{Input} directory. However, any input file can be used by entering the following command in a console:\\ FreeFem++ \textit{BEC\_$\mathcal X$D\_ToolBox.edp} -param \textit{name\_physics} -run \textit{name\_run}.\\ Here, \textit{name\_physics} is the name of the input file containing the physical parameters and \textit{name\_run} is the name of the input file containing the parameters for the computation.\\
1/ 
 The first file in the \textit{Input} directory, \textit{BEC\_$\mathcal X$D\_physic\_param.dat}, contains the physical parameters:
\begin{itemize}
 \item @scaling, a string that can take the values \textit{AR} or \textit{Classical} depending on which scaling is chosen. A value must be given to this parameter.
 \item @kind, a boolean that takes the value 0 if one wants to set constants already built from \eqref{eq-scal-Cg}, \eqref{eq-scal-Com} and \eqref{eq-scal-trap-ax}, or the value 1 if one wants to set the corresponding physical parameters. A value must be given to this parameter.
 \item If the 0 value was chosen for @kind the following parameters must be set to a real value:\\
@beta (=$\beta$), the coefficient in front of the non linear part of the equation (see \ref{eq-scal-Cg}),\\
@Omop ($= \frac{\Omega}{\omegap}$), the coefficient of the angular momentum (see \ref{eq-scal-Com}),\\
@ax, @ay, @az and @a4 are the coefficients in the potential $V_\text{trap}$ (see \ref{eq-scal-trap-ax}).
 \item If the value 1 was given to @kind one must give a real value to the following parameters:\\
 @N, the number of atoms,\\
 @m, the atomic mass,\\
 @as, the scattering length,\\
 @Omega ($=\Omega$), the rotation speed,\\
 @omegax ($=\omega_x$), @omegay ($=\omega_y$), omegaz ($=\omega_z$), @omega2 ($=\omega_2$), @omega4 ($=\omega_4$), @U2 and @U4 are the coefficients in $V_\text{trap}$ (see \ref{eq-GP-trap-V}).
\end{itemize}
2/
 The second file, \textit{BEC\_$\mathcal X$D\_run\_param.dat} contains the parameters for the run:
\begin{enumerate}
 \item Here are the parameters that must be set,
\begin{itemize}
 \item @method is a string to choose a method. The possible values are \textit{Ipopt} or \textit{GradS}.
 \item @EPS0 is a real corresponding to the final error to reach.
 \item @init is a string with the name of the initial approximation to use. The possible values are \textit{TF} (Thomas Fermi), Ipoptaxi (axisymmetric approximation) or  Ipoptnorot (no rotation).
\end{itemize}
\textbf{All the parameters that follow are set by default:}
 \begin{itemize}
 \item @GradSMaxIter is the maximum number of iterations in the Sobolev gradient method. Default value: 8000.
 \item @IpoptMaxIter is the maximum number of iterations between each mesh adaptation in Ipopt method. Default value: 50.
 \end{itemize}
 \item The following parameters are used for the outputs:
\begin{itemize}
 \item @dircase is a prefix of the name of the output directory. The form of potential and the name of the method used for computation are automatically added to this name. Default value: BEC\_3D.
 \item @scase is a prefix of the name of the output files. The values of $\frac{\Omega}{\omegap}$ and $C_g$ are automatically added to this name. Default value: BEC\_3D.
 \item @withplot is a boolean controlling the possibility of plotting the solution during the run. Default value: 1.
 \item @savesol is a boolean controlling the possibility of saving the solution during the run. Default value: 1.
 \item @IWAIT is a boolean controlling the possibility of  waiting after each plot. Default value: 0.
 \item @meditplot is a boolean controlling the possibility of plotting the solution with medit. Default value: 0.
 \item @output is a string that takes the value \textit{vtk} or \textit{tecplot} for the outputs format. Default value: \textit{tecplot}.
 \item @ITERSAVE, @ITERNORM and @ITERPLOT are integers corresponding to the frequency of iterations in Sobolev gradient method to save, normalize or plot the solution. Default value: 100.
 \item @savenergy and @plotenergy are booleans to save and plot the energy during the run. Default value: 1.
 \item @countvortices (only in 2D) is a boolean to count the number of vortices and to give their position. Default value: 1.
\end{itemize}
 \item One can control how to build the initial mesh by setting the following parameters:
\begin{itemize}
 \item @aRdom is a coefficient that multiply the Thomas-Fermi radius in order to have a larger domain. Default value: 1.25.
 \item @nbseg is the number of segments on the border of the mesh. Default value: 50 in 3D and 200 in 2D.
 
 \item @meshkind (only in 3D) is a string that can take the values \textit{cylindre} or \textit{ellipsoid} and allows the user to choose between a cylindrical mesh or an ellipsoidal mesh in 3D. Default value: \textit{ellipsoid}.
 
 \item @hminsurf is the minimal size of the edge of a triangle on the surface of the ellipsoidal mesh. Default value: 0.6. 
 \item @hminvol is the minimal size of the edge of a tetraedra inside the ellipsoidal mesh. Default value: 0.3. 

\end{itemize}
 \item The parameters for loading an old solution as an initial field are:
\begin{itemize}
 \item @ifILrst is a boolean, with true value if the user wants to load a restart file. Default value: 0.
 If ifILrst = 1, the following 4 parameters have to be specified:
 \item @keepmesh is a boolean to choose to keep the loaded mesh or not.
 \item @dirload is a string with the name of the directory containing the restart mesh and solution.
 \item @dmesh is a string with the name of the file containing the mesh to load.
 \item @dsol is a string with the name of the file containing the solution to load.
\end{itemize}
 \item The following parameters control how to build the initial field for the wave function: 
\begin{itemize}
 \item @mod is an integer. If @mod $>0$, a central vortex with winding number @mod is added in the axisymmetric approximation built with Ipopt. Default value: 0.
 \item @narray is the number of circles of vortices in the manufactured initial field (see section \ref{sec-manuf}). Default value: 0.\\
 If narray = 1, the following 8 parameters have to be specified:
 \item @Nv, the number of vortices on each circle.
 \item @Rarr, the radius of the first circle.
 \item @dRarr, the distance between two circles.
 \item @Tharr, the orientation of the first circle.
 \item @dTharr, a step between the orientation of each circle.
 \item @shape (only in 3D) is a string controlling the shape of the ansatz vortex added in the initial condition (see equation \ref{eq-anzatz-shape}). It can be I, Ux, Uy, Sx, Sy. The x or y indicates if the vortex is in the plane (O,x,z) or (O,y,z). Default value: I.
 \item @curvature $(=\alpha_v$, only in 3D) is the parameter controlling the curvature of the ansatz vortex (see \ref{eq-anzatz-shape}). Default value: 10.
 \item @length $(=\beta_v$, only in 3D) is the parameter controlling the length of the ansatz vortex (see \ref{eq-anzatz-shape}).  Default value: 2.
\end{itemize}
 \item The user can control the mesh adaptivity process:
\begin{itemize}
 \item @ifIadapt is a boolean to choose to adapt the mesh of the initial field. Default value: 1. 
 \item @erradaptI is the error in the mesh adaptation of the initial field. This parameter is used by the FreeFem function adaptmesh in 2D or mshmet in 3D. Default value: 0.01 in 3D and 0.1 in 2D.
 \item @ifRadapt is a boolean to choose to adapt the mesh  during the computation. Default value: 1.
 \item @hminad is the minimal size of an edge in the new mesh. Default value: 0.001.
 \item @hmaxad is the maximal size of an edge in the new mesh. Default value: 1.
 \item @erradapt is the error in the mesh adaptation. It does change during a computation with Sobolev gradient method. If the Ipopt method is used for the computation, it corresponds to the parameter $\varepsilon_\text{last}$ in (\ref{eq-algo-ipopt}). Default value: 0.01 in 3D and 0.1 in 2D for Sobolev gradient method and 0.008 in 3D and 0.005 in 2D for Ipopt method.
 \item @anisoadapt is a real value. If @anisoadapt $ > 0$, the mesh adaptation will be anisotropic and the ration between the size of the smallest and the biggest edges of each triangle will be bounded by @anisoadapt. Default value: 10.
\end{itemize}
 \item The following parameters are needed for the mesh adaptation in the Sobolev gradient method only:
\begin{itemize}
 \item @EPSAD1 is the first value of the $L^2$ relative error the user wants to reach to make a mesh adaptation (= $\varepsilon^1$ in \ref{sec-mesh-adapt} ). Default value: 1e-2.
 \item @EPSADMIN is the last stage (= $\varepsilon_c$ in \ref{sec-mesh-adapt} ). Default value: 1e-9.
 \item @IPASSAL is the number of times a mesh adaptation is performed before changing the value of EPSAD1 ($N_{ad}$ in \ref{sec-mesh-adapt} ). Default value: 2.
 \item @EPSADSTEP is a factor to change the value of EPSAD1. Default value: 2.
 \item @ITERADAPT is the maximum number of iterations between two mesh adaptations. If $ITERADAPT= 0$, we don't use this criterion. Default value: 0.
\end{itemize}
 \item The last parameters are for the mesh adaptation in Ipopt method:
\begin{itemize}
 \item @niadapt is the number of times a mesh adaptation is performed with the same error $\varepsilon_k$ (see \ref{eq-algo-ipopt}). Default value: 1.
 \item @nbadapt is the total number of mesh adaptations made during the computation ($n_\text{adapt}$ in \ref{eq-algo-ipopt}). Default value: 6 in 3D and 4 in 2D.
 \item @maerr1 is the initial error in mesh adaptation ($\varepsilon_0$ in \ref{eq-algo-ipopt}). Default value: 0.01.
\end{itemize}
\end{enumerate}

\subsection{Output files}
When a computation starts, the \textit{Output} directory is created.
It contains a directory, whose name includes the prefix (defined by the parameter @dircase), the form of potential and the chosen method. 
This directory will contain an \textit{.echo} file with a summary of the main parameters, informations on the run, names of the output files, final energy and the CPU time.
The \textit{plot.gp} file will contain a Gnuplot script that the user can run to plot the evolution of the energy, the error, the angular momentum or the $L^2$ norm of the solution.
The \textit{.mesh} and \textit{.rst} file contains the mesh and the solution respectively. They can be used as a restart field.
Finally, the \textit{.tec} or \textit{.vtk} files contain the solution for a given iteration (defined by the parameters ITERPLOT) in the format \textit{tecplot} or \textit{vtk}.

\begin{figure}[ht]
\begin{center}
\includegraphics[width=0.9\textwidth]{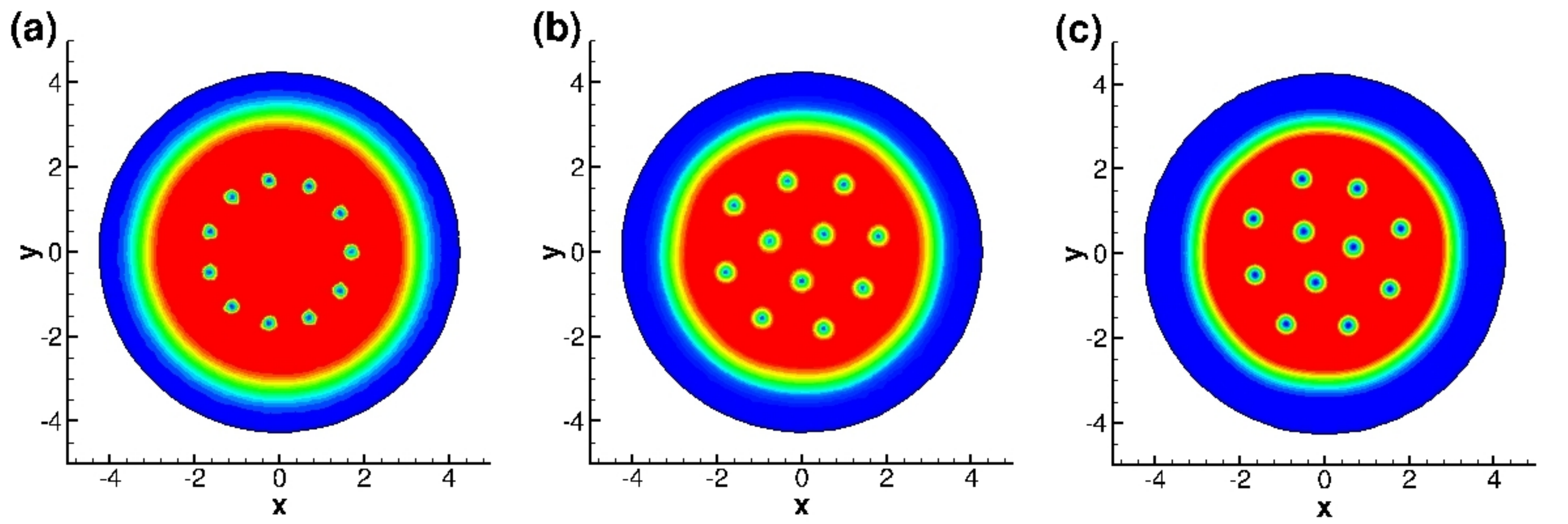}
\end{center}
\caption{Solution of the second example in 2D: (a) the initial state built by adding 11 manufactured vortices to the ground state computed with Ipopt (1D axisymmetric), (b) converged solution with the Sobolev gradient method, (c) solution obtained with Ipopt. }
\label{fig-2D-example}
\end{figure}

\section{Examples and user interface}\label{sec-example}
To simplify the understanding of parameter files, some examples are provided in the directory \textit{Examples}. A user interface was implemented using the \textbf{GLUT} library \cite{glut} to run these examples or to run the toolbox with predefined parameter files. In this section, we first present the examples files and some results of computations. Then, we focus on the use of the GLUT user interface. 
\subsection{2D computations}
The examples for 2D computations use two forms of the trapping potential. For each case, the use of both scalings and numerical methods is possible.
\begin{enumerate}
\item The first case is the harmonic potential with $a_x = 1, \, a_y = 1$ (see \ref{eq-scal-trap-ax}), $\beta = 500$ (see eq. \ref{eq-scal-Cg}) and $\Omega/\omegap = 0.4$ (see eq. \ref{eq-scal-Com}).
We start with an initial approximation made with Ipopt axisymmetric and we add one manufactured off-centred vortex, as in Figure \ref{fig-2D-adapt}(a).
The final state we reach is a BEC with one central vortex as in figure \ref{fig-2D-adapt}(b). To run this example, the following files from the directory \textit{Examples} have to be used:
\begin{itemize}
\item \textit{AR\_Harm\_physic\_param.dat} or 
\textit{Classical\_Harm\_physic\_param.dat} for the physical parameters, depending on which scaling is chosen,
\item \textit{Ipopt\_Harm\_run\_param.dat} or \textit{GradS\_Harm\_run\_param.dat} for the computation parameters, depending on which method is chosen.
\end{itemize} 
\item The second case is a combined quartic/quadratic potential with $a_x = 1, \, a_y = 1, \, a_4 = 0.5$ (see eq. \ref{eq-scal-trap-ax}), $\beta = 500$ and $\Omega/\omegap = 2$.
We start with an initial approximation made with Ipopt axisymmetric and we add a circle of manufactured vortices, as in figure \ref{fig-2D-example}(a). 
Both methods reach a BEC with eleven vortices organised into an Abrikosov lattice as shown in figures \ref{fig-2D-example}(b) and \ref{fig-2D-example}(c). To run this example, the following files must be selected by the user in the directory \textit{Examples}:
\begin{itemize}
\item \textit{AR\_Quart\_physic\_param.dat}
 or \textit{Classical\_Quart\_physic\_param.dat} for the physical parameters,
\item \textit{Ipopt\_Quart\_run\_param.dat} 
 or \textit{GradS\_Quart\_run\_param.dat} for the chosen method.
 \end{itemize}
\end{enumerate}

\begin{figure}[h]
\begin{center}
\includegraphics[width=\textwidth]{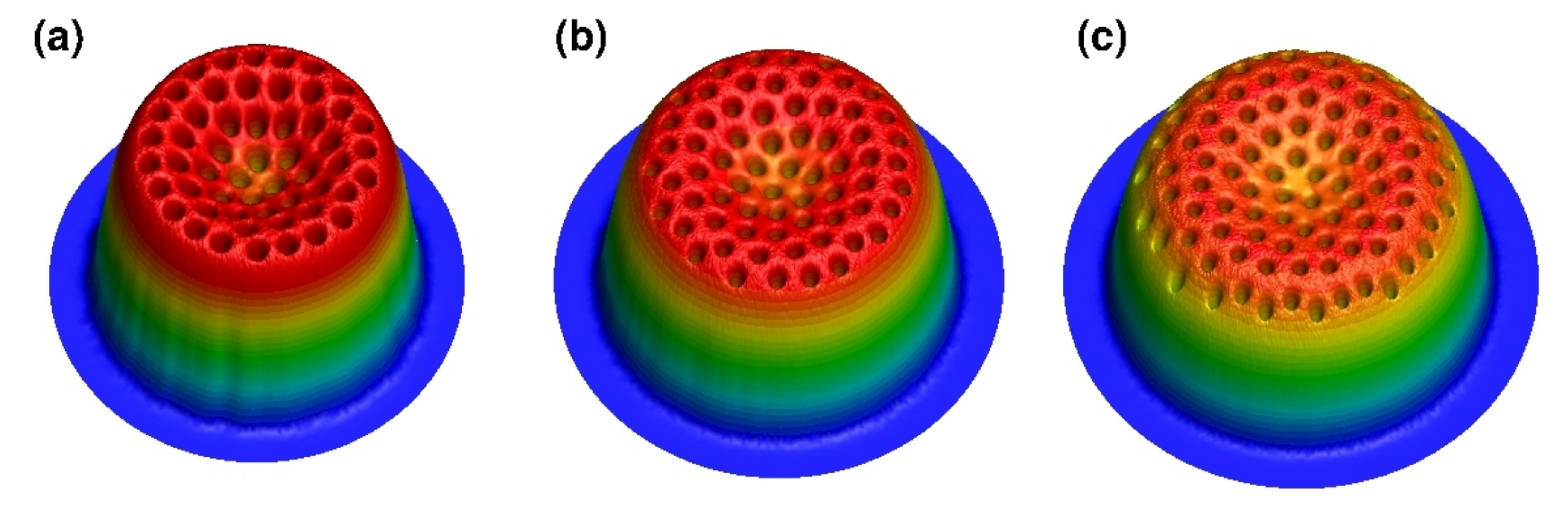}
\end{center}
\caption{2D solution obtained with the Sobolev gradient method for a quartic+quadratic potential with $\Omega/\omegap = 3.5$ and different values of the non-linear interaction constant: $(a)\,  \beta=5000,\,(b)\, \beta=10000,\,(c)\, \beta=15000$.}
\label{fig2DLat}
\end{figure}

\begin{figure}[h]
\begin{center}
\includegraphics[width=\textwidth]{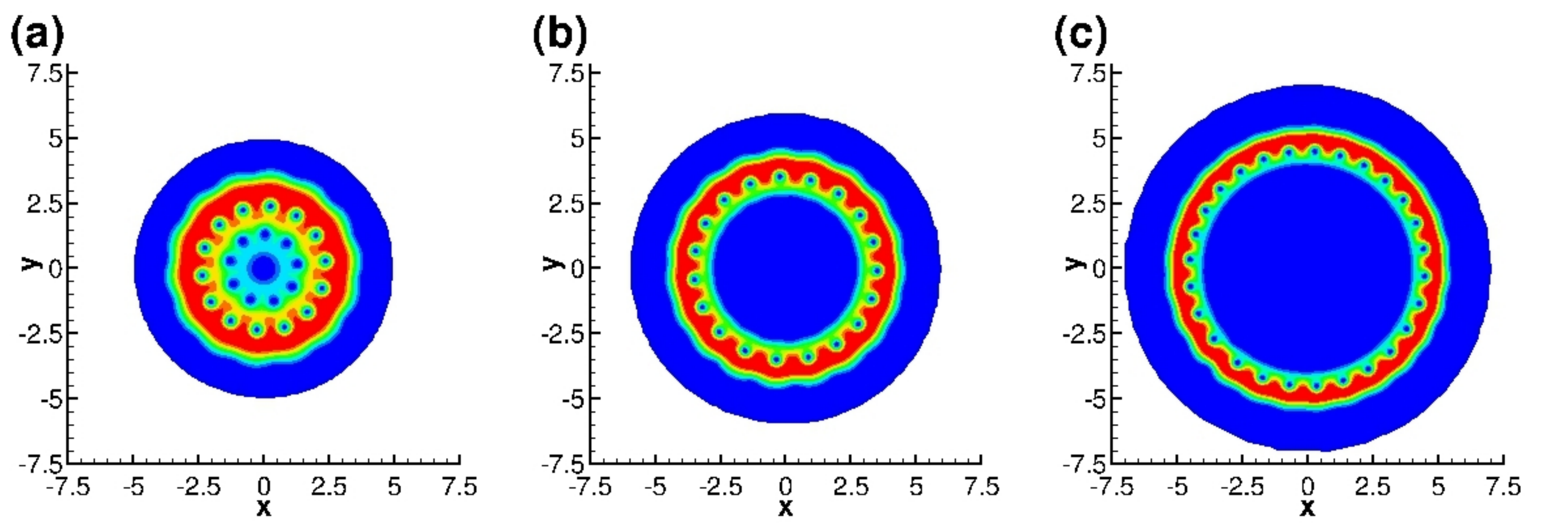}
\end{center}
\caption{2D solution built with the Ipopt method for a quartic+quadratic potential with $\beta = 500$ and different values of the rotation frequency: $(a) \, \Omega/\omegap=3,\, (b)\, \Omega/\omegap=4,\, (c) \,  \Omega/\omegap=5$.}
\label{fig2DGiant}
\end{figure}

In figures \ref{fig2DLat} and \ref{fig2DGiant}, we provide two results with the same quartic+quadratic potential as in the 2D example illustrated in figure \ref{fig-2D-example}: $a_x = 1,\, a_y = 1,\, a_4 = 0.5$. In the case illustrated in figure \ref{fig2DLat}, we set the rotation speed to $\Omega = 3.5$ and increase the non-linear constant $\beta$ from $5000$ to $15000$. When this constant increases, the condensate becomes larger and  the number of vortices increases significantly. They arrange in a triangular Abrikosov lattice. The files used to perform this simulation are provided in the directory \textit{Input} as \textit{BEC\_2D\_physic\_param\_Latt.dat} and \textit{BEC\_2D\_run\_param\_Latt.dat}.

In the case of figure \ref{fig2DGiant}, the non-linear constant $\beta = 500$ is fixed and the rotation speed $\Omega$ increases from $3$ to $5$. The condensate is larger when the rotation speed increases and a giant vortex appears at the centre of the condensate. This case was simulated in \cite{BEC-physV-2005-fetter}. The size of the computational domain increases as the rotation speed increases. This illustrates the need of the use of the Thomas-Fermi approximation to estimate the size of the domain. 
The files used to perform this simulation are provided in the directory \textit{Input} under the names \textit{BEC\_2D\_physic\_param\_Giant.dat} and \textit{BEC\_2D\_run\_param\_Giant.dat}.

\begin{figure}[h]
\begin{center}
\includegraphics[width=0.9\textwidth]{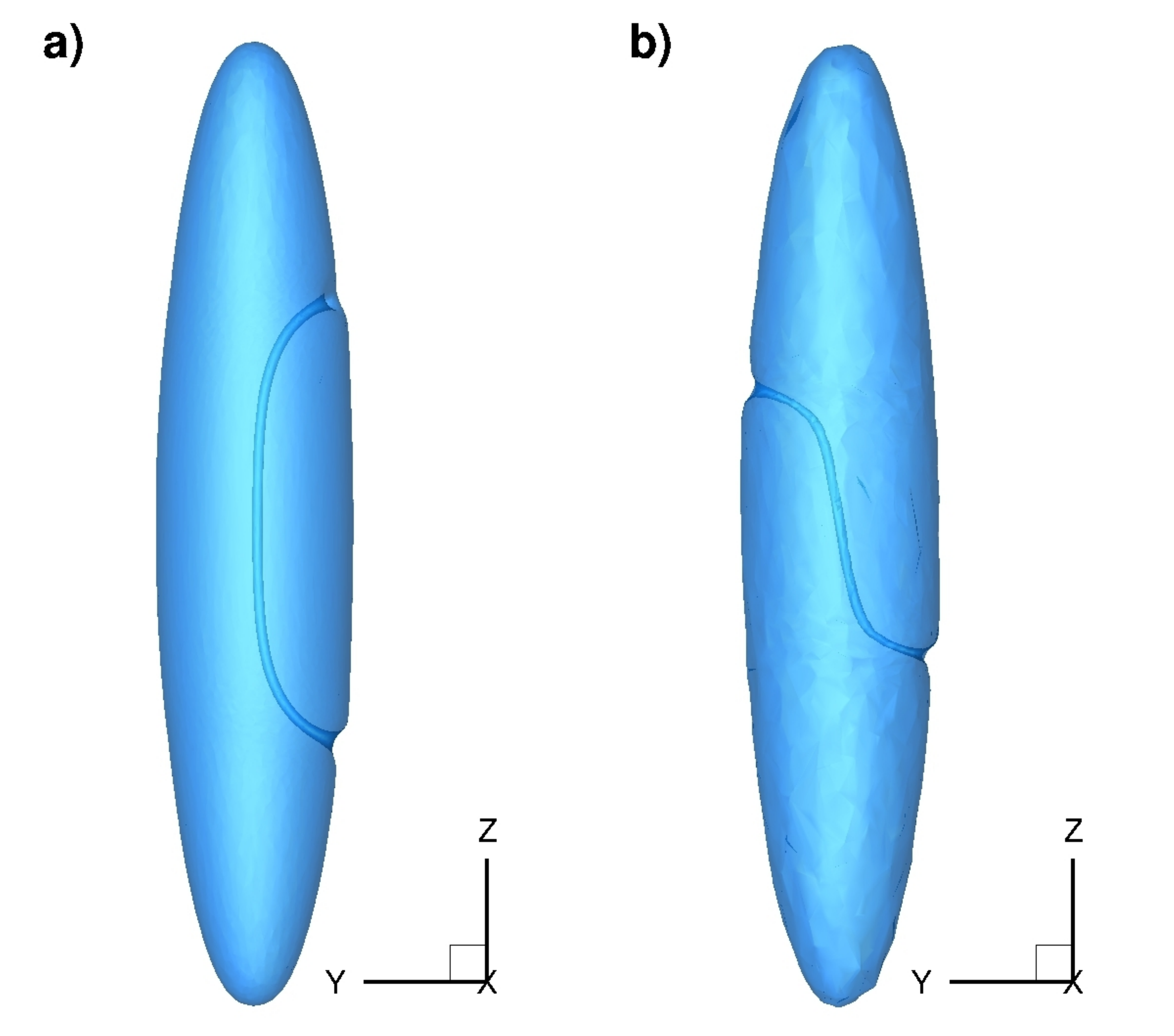}
\end{center}
\caption{Example from the interface software: isosurface of low atomic density illustrating 3D vortices. (a) U vortex obtained with the Ipopt method starting from a manufactured initial state with a U-shaped centred vortex. (b) S vortex obtained with the Sobolev gradient method starting from a manufactured initial state with a S-shaped centred vortex. }
\label{fig-3D-example}
\end{figure}

\subsection{3D computations}

Two examples with harmonic trapping potential with $a_x = 1, \, a_y = 1.06^2,\, a_z = 0.067^2$ (see eq. \ref{eq-scal-trap-ax}), $\beta = 15900$ and $\Omega/\omegap = 0.4$ are provided for 3D computations. They correspond to numerical tests used in \cite{dan-2003-aft}. These tests have shown that, with the same physical parameters, one can get different final meta-stable states, by starting from different initial states. In the first case, illustrated in figure \ref{fig-3D-example}(a), the computation starts with an axisymmetric approximation with a manufactured U-shaped vortex added at the centre. The final state, reached with both the Sobolev gradient method and the Ipopt method, presents a bended vortex with a U shape. In the second case of figure \ref{fig-3D-example}(b), we start with an axisymmetric approximation with a manufactured central vortex with a S shape. The final converged state keeps a S-shaped vortex when using both numerical methods. According to \cite{dan-2003-aft} the S-shaped vortices is a local minima of the energy. We conclude that both methods converge to the local minimum which is the closest to the initial guess provided. The input files used for these examples are provided in the directory \textit{Examples} as:
\begin{itemize}
  \item \textit{VortexU\_physic\_param.dat} and \textit{VortexS\_physic\_param.dat} for the physical parameters,
\item  \textit{VortexU\_GradS\_run\_param.dat} and  \textit{VortexU\_Ipopt\_run\_param.dat} for the computation parameters of the U-shaped vortex case,
\item \textit{VortexS\_GradS\_run\_param.dat} and \textit{VortexS\_Ipopt\_run\_param.dat} for the computation parameters of the S-shaped vortex case.
\end{itemize}  

\begin{figure}[h]
\begin{center}
\includegraphics[width=\textwidth]{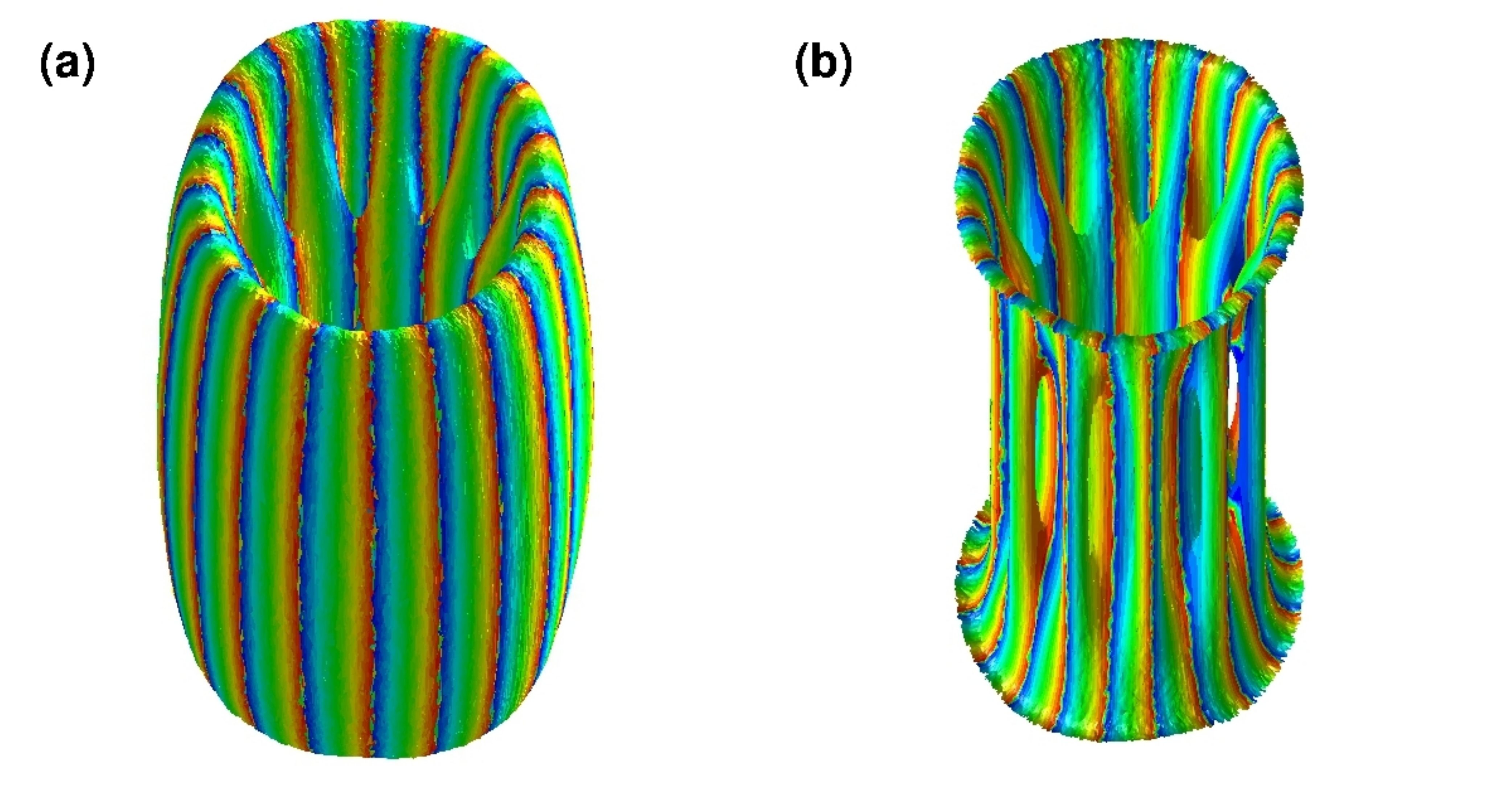}
\end{center}
\caption{3D solution computed with the Ipopt method for a quartic-minus-quadratic potential.  (a) Isosurface of low atomic density coloured with the phase. (b) The outer layer is removed to see the singly quantised vortices inside the condensate. $\beta = 21000$ and $\Omega/\omegap = 2$.}
\label{fig3DGiant}
\end{figure}

The result shown in figure \ref{fig3DGiant} was obtained using physical parameters from \cite{dan-2004-aft}: $a_x = -0.2, \, a_y = -0.2,\, a_z = 0.067^2, \, a_4 = 0.075 , \, \beta = 21000$ (see \ref{eq-scal-trap-ax}), and $\Omega/\omegap = 2$. The ground state displays a giant vortex surrounded by eleven singly-quantized vortices. This simulation was carried out using Ipopt for a quartic-minus-quadratic potential. The files used to perform this simulation are provided in the directory \textit{Input} as: \textit{BEC\_3D\_physic\_param\_Giant.dat} and\\ \textit{BEC\_3D\_run\_param\_Ipopt\_Giant.dat}.

\begin{figure}[h]
\begin{center}
\includegraphics[width=0.9\textwidth]{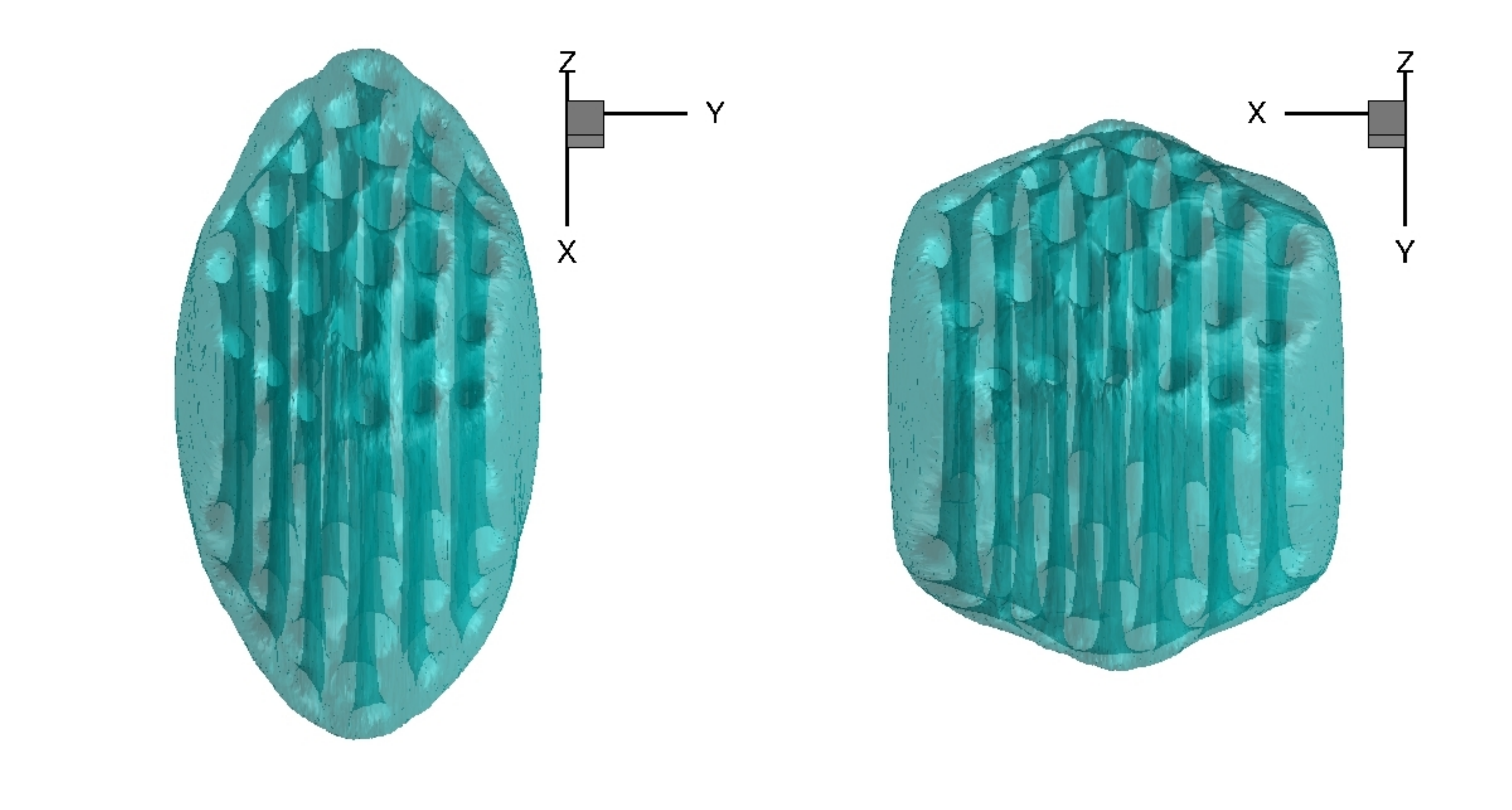}
\end{center}
\caption{3D solution computed with the Sobolev gradient method for an anisotropic harmonic potential.  Different views of an isosurface of the low atomic density showing the presence of 31 singly quantized vortices in an elongated condensate. $\beta = 50000$ and $\Omega/\omegap = 0.95$.}
\label{fig3DLattice}
\end{figure}

Figure \ref{fig3DLattice} illustrates other possible vortex states that can be obtained. An anisotropic harmonic potential with $a_x = 1, \, a_y = 1.06^2,\, a_z = 0.067^2, \, \beta = 50000$ (see \ref{eq-scal-trap-ax}), and $\Omega/\omegap = 0.95$  was used. The computation resulted in an Abrikosov lattice with 31 vortices in 3D. The anisotropy makes the condensate to take an elongated shape following the x-axis. The files used to perform this simulation are provided in the directory \textit{Input} as: \textit{BEC\_3D\_physic\_param\_aniso.dat} and \textit{BEC\_3D\_run\_param\_Ipopt\_aniso.dat}.

\subsection{Optional user interface}
\begin{figure}[h]
\begin{center}
\includegraphics[width=0.9\textwidth]{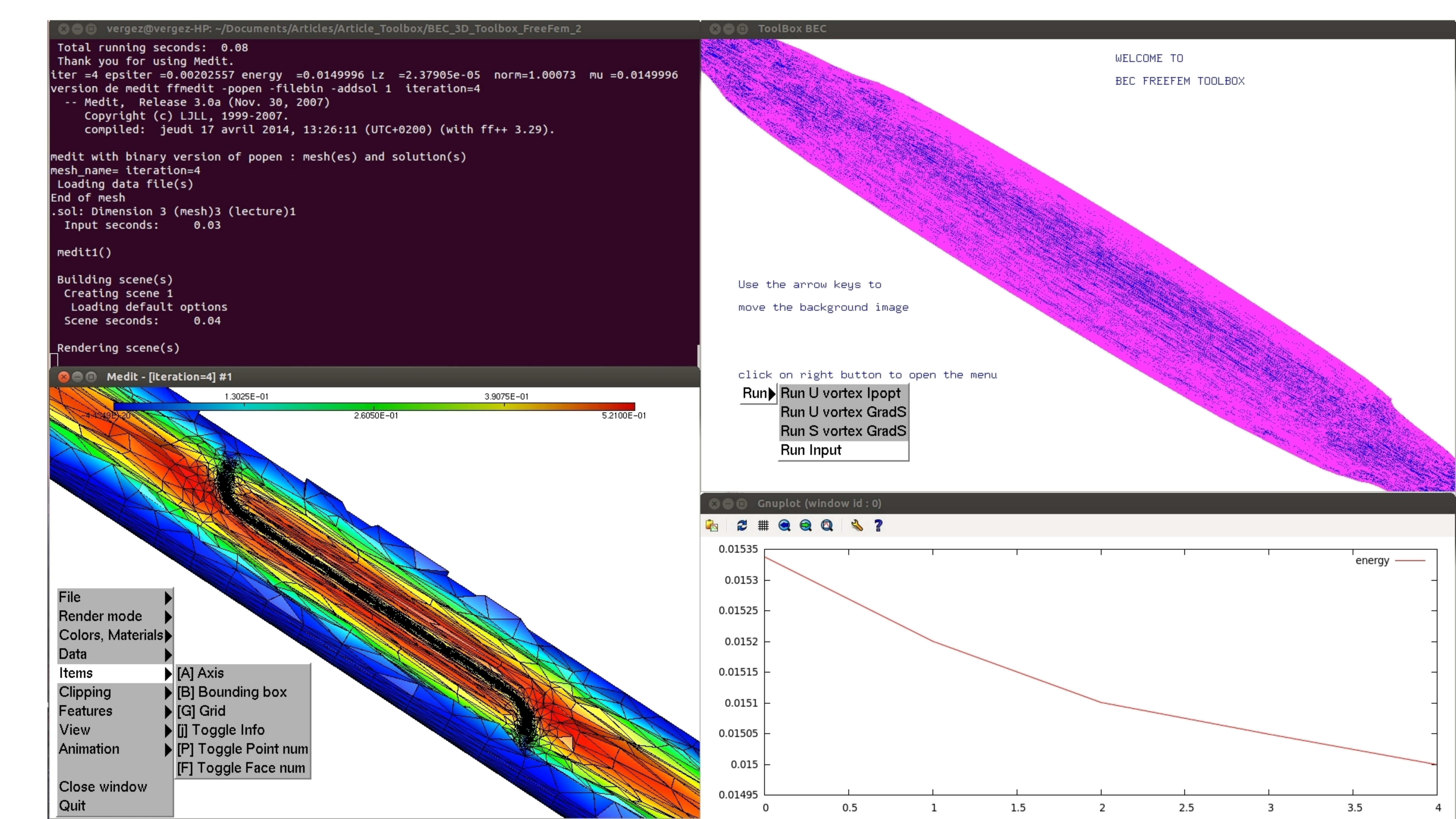}
\end{center}
\caption{Screen capture of the user interface for a 3D computation. We can see: the terminal in which the application is run, the solution plotted with Medit and the plot of the energy evolution.}
\label{figInter}
\end{figure}

A simple user interface was made in C++ with the GLUT tool of the OpenGL library. The C++ source code for this interface is in the directory GLUT. It can be compiled using the makefile provided with the toolbox. This interface allows the user to easily run  the examples or to run the toolbox using any modified input file.
The screen capture of the interface in 3D is shown in figure \ref{figInter}. On the top left corner one can see a terminal from which was run the executable "RunToolbox". The window on the top right corner of figure \ref{figInter} appears. By clicking on the right button of the mouse, a pull-down menu allows the user to run the toolbox with one of the three example files provided, or using the input files from the \textit{Input} directory.
Then a Gnuplot window appears plotting the evolution of the energy during the run. This window is on the bottom right corner of figure \ref{figInter}. Finally the bottom left corner of figure \ref{figInter} shows the 3D solution plotted with Medit. The user can also decide to plot it with the usual graphical interface of FreeFem++.

In 2D, there are three menus to choose one of the examples previously described:
\begin{enumerate}
\item \textit{Potential} allows to choose between the harmonic (\textit{Harm}) or the quartic+quadratic (\textit{Quart}) trapping potential example.
\item \textit{Method} allows to choose between the Ipopt (\textit{Ipopt}) or the Sobolev gradient (\textit{GradS}) method.
\item \textit{Scaling} allows to choose between the Aftalion Riviere (\textit{AR}) or the classical scaling.
\end{enumerate}
In the last menu (\textit{Run}), the user can run either the selected example (\textit{Run Example}) or run the toolbox with the input files which are in the \textit{Input} directory (\textit{Run Input}).\\

\section{Conclusion}
We provide with this paper a finite-element software for 2D and 3D computation of stationary solutions of the Gross-Pitaevskii equation. The user has the choice between two robust and optimised numerical methods: a steepest descent method based on Sobolev gradients and a minimization algorithm based on the state-of-the-art optimization library Ipopt. For both methods, mesh adaptivity strategies are implemented to reduce the computational time and increase the local spatial accuracy when vortices are present. The numerical system is tested and validated through various cases representing 2D and 3D configurations of Bose-Einstein condensates in rotation. A particular attention was paid to the physical interpretation of the computations. 
The main parameters of the run can be prescribed either in non-dimensional or physical form. Thomas-Fermi approximations are derived as closed formulae for a more accurate description of the initial field for the minimization procedures.  Energy and angular momentum are tracked during the computation and post-processing tools allow to identify quantized vortices in the final, equilibrium state.  

An optional graphical user interface is also provided with the software. It allows to easily run predefined cases or with user-defined parameter files.  

The programs were written as a toolbox to be used within  the free software FreeFem++. This offers the advantage that all technical issues related to the implementation of the finite element method are hidden, allowing to focus on numerical algorithms and their performance. Automatic mesh generators, powerful mesh adaptivity functions and the  availability of various types of finite elements with complex functions are the main features making FreeFem++ very appealing in implementing numerical methods for Schr\"odinger type equations.  The toolbox distributed with this paper is extremely versatile and can be easily adapted to deal with different physical models. A natural extension of this toolbox is the simulation of the time-dependent Gross-Pitaevskii equation: this is an ongoing work and will be reported in a further contribution.


\section*{Acknowledgements}
This work was supported by the French ANR grant  ANR-12-MONU-0007-01 BECASIM (Mod{\'e}les Num{\'e}riques call).
We would like to acknowledge the use of computational resources provided by CRIHAN (Centre de Ressources Informatiques de Haute-Normandie, France) under the project 2015001.

\appendix

\section{Formulae for the Thomas-Fermi approximation}
We derive in this Appendix closed relationships for the Thomas-Fermi approximation for different types of trapping potentials (quartic $\pm$ quadratic). The Thomas-Fermi density (\ref{eq-scal-TF}) can be rewritten using (\ref{eq-scal-Ctrap-eff}) as:
\begin{equation}
\label{eq-app-TF}
\rtf = \mod{u}^2= \left(\frac{\rho_0 - 2\vadim{V}^\text{eff}}{C_\text{S}}\right)_+, \quad \rho_0 =  2\varepsilon  {\vadim{\mu}}  = 2\varepsilon \frac{\mu}{\hbar \omegap}, \quad \mbox{and} \quad C_\text{S} = 2\varepsilon^2C_g.
\end{equation} 
We recall that $\varepsilon=1$ for the classical scaling using the oscillator length $\aho$ as length scale. The constant $\rho_0$ will be determined by imposing the unitary norm constraint (\ref{eq-scal-cons}). We derive below different formulas for $\rho_0$ corresponding to the effective trapping potential (\ref{eq-scal-trap-V}). We drop in the following the tilde notation.

\subsection{2D harmonic potential}
For this case, the effective trapping potential (\ref{eq-scal-trap-V}) is reduced to
\begin{equation}
\vadim{V}^\text{eff} = \frac12\left(a_x x^2 + a_y y^2\right).
\end{equation}
The unitary norm constraint (\ref{eq-scal-cons}) becomes
\begin{equation}
\label{eq-app-calc-rho2}
I=\int_{\mathcal D} \left(\rho_{0} - a_x x^2 - a_y y^2\right)\, dxdy = C_\text{S}.
\end{equation}
To calculate $I$ analytically, we use the change of variables:
\begin{equation}
\left\{
\begin{array}{lcl}\vspace{0.2cm}
x &=& \displaystyle \frac{r}{\sqrt{a_x}}\, \cos \theta\\
y &=& \displaystyle \frac{r}{\sqrt{a_y}}\, \sin \theta
\end{array}
\right.
\quad
dxdy = \frac{r}{\sqrt{a_x a_y}}\, drd\theta, \quad r\in [0, \sqrt{\rho_{0}}], \,\, \theta\in [0, 2\pi],
\end{equation}
and
\begin{equation}
I=\frac{1}{\sqrt{a_x a_y}} \int_0^{2\pi} d\theta \,\int_0^{\sqrt{\rho_{0}}}\, (\rho_{0}-r^2)r dr = \frac{\pi \rho_{0}^2}{2 \sqrt{a_x a_y}}.
\end{equation}
Finally, the constant $\rho_{0}$ is expressed as:
\begin{equation}
\label{eq-app-h2D-rho0}
\rho_{0} = \left( \frac{2 \sqrt{a_x a_y}}{\pi}\, C_\text{S}\right)^{1/2},
\end{equation}
and the dimensions of the condensate follow:
\begin{equation}
\label{eq-app-h2D-Rxy}
R_x = \sqrt{\frac{\rho_{0}}{a_x}}, \quad R_y = \sqrt{\frac{\rho_{0}}{a_y}}.
\end{equation}

\subsection{3D harmonic potential}

Same analysis for the potential
\begin{equation}
\vadim{V}^\text{eff} =  \frac12\left(a_x x^2 + a_y y^2 + a_z z^2\right).
\end{equation}
The constraint \eqref{eq-scal-cons} becomes
\begin{equation}
\label{eq-app-calc-rho3}
I=\int_{\mathcal D} \left(\rho_{0} - a_x x^2 - a_y y^2 - a_z z^2 \right)\, dxdydz = C_\text{S}.
\end{equation}
To calculate $I$ analytically, wa use the change of variables:
\begin{equation}
\left\{
\begin{array}{lcl}\vspace{0.2cm}
x &=& \displaystyle \frac{r}{\sqrt{a_x}}\, \sin\theta \cos\phi\\ \vspace{0.2cm}
y &=& \displaystyle \frac{r}{\sqrt{a_y}}\, \sin\theta \sin\phi\\
z &=& \displaystyle \frac{r}{\sqrt{a_z}}\, \cos\theta
\end{array}
\right.
\quad
dxdydz = \frac{r^2 \sin\theta}{\sqrt{a_x a_y a_z}}\, drd\theta d\phi, \quad
\left\{
\begin{array}{lcl}\vspace{0.2cm}
r &\in& \displaystyle [0, \sqrt{\rho_{0}}]\\ \vspace{0.2cm}
\theta &\in& \displaystyle [0, \pi]\\
\phi &\in& \displaystyle [0, 2\pi]
\end{array}
\right.
\end{equation}
and
\begin{equation}
I=\frac{1}{\sqrt{a_x a_y a_z}} \int_0^{2\pi} d\phi  \int_0^{\pi}  \sin\theta d\theta\,\int_0^{\sqrt{\rho_{0}}}\, (\rho_{0}-r^2)r^2 dr = \frac{8 \pi \rho_{0}^{5/2}}{15 \sqrt{a_x a_y a_z}}.
\end{equation}
Finally, the constant $\rho_{0}$ is expressed as:
\begin{equation}
\label{eq-app-h2D-rho03D}
\rho_{0} = \left( \frac{15 \sqrt{a_x a_y a_z}}{8\pi}\, C_\text{S}\right)^{2/5},
\end{equation}
and the dimensions of the condensate follow:
\begin{equation}
\label{eq-app-h2D-Rxy3D}
R_x = \sqrt{\frac{\rho_{0}}{a_x}}, \quad R_y = \sqrt{\frac{\rho_{0}}{a_y}},
\quad R_z = \sqrt{\frac{\rho_{0}}{a_z}}.
\end{equation}

\subsection{2D combined quartic and quadratic potential}

We consider that the trap has radial symmetry ($a_x=a_y=a_2$) and the trapping potential is
\begin{equation}
\vadim{V}^\text{eff} =  \frac12\left(a_2 r^2 + a_4 r^4\right).
\end{equation}
Note that $a_4>0$, but $a_2$ can be either positive (quartic+quadratic potential) or negative (quartic-quadratic potential).
The border of the condensate is defined by the radius $R$ that satisfies:
\begin{equation}
\label{eq-app-q2D-roots}
a_4 R^4 + a_2 R^2 -\rho_{0} =0, \Longrightarrow R^2_{\pm} = \frac{-a_2 \pm \sqrt{a_2^2+4\rho_{0} a_4}}{2 a_4}.
\end{equation}

\subsubsection{Case $a_2 \ge 0$: quartic+quadratic potential}

In this case, $a_2>0, a_4>0$ and we infer from \eqref{eq-app-q2D-roots} that $\rho_{0}>0$ and it exists a single root $R_+$:
\begin{equation}
R^2_+ = \frac{-a_2 + \sqrt{a_2^2+4\rho_{0} a_4}}{2 a_4} > 0.
\end{equation}

The constraint \eqref{eq-scal-cons} becomes in polar coordinates $(r,t)$:
\begin{equation}
I=2\pi \int_{0}^R \left(\rho_{0} - a_2 r^2 - a_4 r^4\right)\, rdr = C_\text{S},
\end{equation}
or
\begin{equation}
C_\text{S} = 2\pi \left[\rho_{0} \frac{R^2}{2}- a_2 \frac{R^4}{4} - a_4 \frac{R^6}{6}\right] = \frac{\pi R^4}{6}\left(3 a_2 + 4 a_4 R^2\right).
\end{equation}
Consequently, we first have to  calculate the root $\eta=R^2>0$ of the non-linear equation (by a Newton method by example):
\begin{equation}
4 a_4 \eta^3 + 3 a_2 \eta^2 - \underbrace{\frac{6}{\pi} \left(C_\text{S} \right)}_{A_\eta} = 0,
\end{equation}
and then calculate
\begin{equation}
\rho_{0} = a_2 \eta + a_4 \eta^2 > 0.
\end{equation}
The radius of the condensate is finally given by $R=\sqrt{\eta}$.

\subsubsection{Case $a_2<0$: quartic-quadratic potential}

We distinguish two cases:
\begin{itemize}
  \item if $\rho_{0}>0$, \eqref{eq-app-q2D-roots} has only one root $R_+$, and the computation of $\rho_{0}$ is the same as above. We also infer that this case occurs when:
      \begin{equation}
      0<\rho_{0} = a_2 R^2 + a_4 R^4 \Longrightarrow \eta = R^2 > \frac{|a_2|}{a_4}.
      \end{equation}

  \item if $\rho_{0}<0$, \eqref{eq-app-q2D-roots} has two roots $R_-, R_+$ and the integration will be carried for $r\in [R_-, R_+]$ (there is a hole in the center of the condensate):
      \begin{equation}
I=2\pi \int_{R_-}^{R_+} \left(\rho_{0} - a_2 r^2 - a_4 r^4\right)\, rdr = C_\text{S},
\end{equation}
or
\begin{equation}
\frac{C_\text{S}}{2\pi}=
 \left[ \frac{\rho_{0}}{2}\left(R_+^2-R_-^2\right)
 -\frac{a_2}{4}\left(R_+^4-R_-^4\right)
 -\frac{a_4}{6}\left(R_+^6-R_-^6\right)\right]
\end{equation}
and using
\begin{equation}
R_+^2 + R_-^2 = -\frac{a_2}{a_4},\,\, R_+^2 R_-^2 =- \frac{\rho_{0}}{a_4},\,\,
R_+^2 - R_-^2 = \frac{\sqrt{a_2^2+4\rho_{0} a_4}}{a_4}
\end{equation}
we obtain that
\begin{equation}
\frac{C_\text{S}}{2\pi}= \frac{\left(a_2^2+4\rho_{0} a_4\right)^{3/2}}{12 a_4^2},
\end{equation}
and finally
\begin{equation}
\rho_{0} = \frac{1}{4 a_4}\left[\left(\frac{6 a_4^2}{\pi} C_\text{S}\right)^{2/3} - a_2^2\right] = \frac{1}{4 a_4}\left[\left(a_4^2 A_\eta\right)^{2/3} - a_2^2\right].
\end{equation}
Since $\rho_{0}<0$, this occurs if
\begin{equation}
a_4 < \sqrt{\frac{\pi |a_2|^3}{6 C_\text{S}}} = \frac{|a_2|^{3/2}}{\sqrt{A_\eta}}.
\end{equation}
\end{itemize}

\subsubsection{Summary for the 2D combined quartic and quadratic potential}

\begin{equation}
\vadim{V}^\text{eff} = \frac12\left(a_2 r^2 + a_4 r^4\right).
\end{equation}
\begin{equation}
\begin{array}{l}
--> \mbox{compute $\rho_{0}$}\\ \vspace{0.2cm}
\left[
\begin{array}{l}\vspace{0.2cm}
\text{Compute} \quad C_\text{S} = 2\varepsilon^2C_g,\\
\mbox{$\bullet$ if $a_2 <0$ and $\displaystyle a_4 < \sqrt{\frac{\pi a_2^3}{6 C_\text{S}}}$} \\ \vspace{0.2cm}
\displaystyle \rho_{0} = \frac{1}{4 a_4}\left[\left(\frac{6 a_4^2}{\pi} C_\text{S}\right)^{2/3} - a_2^2\right] \\ \vspace{0.2cm}
\mbox{$\bullet$ else}\\  \vspace{0.2cm}
\mbox{calculate the root $\eta>0$ of:}\\ \vspace{0.2cm}
\displaystyle f(\eta)=4 a_4 \eta^3 + 3 a_2 \eta^2 - \frac{6}{\pi} \left(C_\text{S} \right) = 0,\\ \vspace{0.2cm}
\mbox{and then calculate}\\
\displaystyle \rho_{0} = a_2 R^2 + a_4 R^4
\end{array}
\right.\\
--> \mbox{compute the maximum radius of the condensate}\\ \vspace{0.2cm}
\displaystyle R_+ = \left(\frac{-a_2 + \sqrt{a_2^2+4\rho_{0} a_4}}{2 a_4}\right)^{1/2}.
\end{array}
\end{equation}

\subsection{3D combined quartic and quadratic potential}

We consider a trapping potential with radial symmetry ($a_x=a_y=a_2$):
\begin{equation}
\vadim{V}^\text{eff} = \frac12\left(a_2 r^2 + a_4 r^4 + a_z z^2\right).
\end{equation}
Note that $a_4>0, a_z >0$, but $a_2$ can be either positive (quartic+quadratic potential) or negative (quartic-quadratic potential).
The border of the condensate is defined by:
\begin{equation}
\label{eq-app-q3D-roots}
a_4 R^4 + a_2 R^2 + a_z z^2 -\rho_{0} =0, \Longrightarrow z(r)= \pm \frac{1}{\sqrt{a_z}}\, \left(\rho_{0} -a_2 r^2 -a_4 r^4\right)^{1/2}.
\end{equation}

\subsubsection{Case $a_2>0$: quartic+quadratic potential}

In this case, $a_2>0, a_z>0, a_4>0$ and we infer from \eqref{eq-app-q3D-roots} that $\rho_{0}>0$ and it exists a single root $R_\perp$ which is the radius of the condensate in the central plane ($z=0$):
\begin{equation}
R^2_\perp = \frac{-a_2 + \sqrt{a_2^2+4\rho_{0} a_4}}{2 a_4} > 0.
\end{equation}
Consequently, the condensate extends in the central plane from $r=0$ to $r=R_\perp$. Using the $z$--symmetry of the condensate, we calculate in cylindrical coordinates
\begin{equation}
I=\int_{{\mathcal D}} (\rho_{0}-a_2 r^2-a_4 r^4-a_z z^2) = \int_0^{2\pi}d\theta\, \int_0^{R_\perp} rdr \, 2\int_0^{z(r)} (\rho_{0}-a_2 r^2-a_4 r^4-a_z z^2) dz
\end{equation}
or using  \eqref{eq-app-q3D-roots}:
\begin{eqnarray} \vspace{0.2cm}
I &=&\displaystyle \frac{8\pi}{3\sqrt{a_z}} \int_0^{R_\perp} (\rho_{0}-a_2 r^2-a_4 r^4)^{3/2} rdr\\
&=&\displaystyle \frac{8\pi}{3\sqrt{a_z}} \int_0^{R_\perp}
\left[\left(\rho_{0}+\frac{a_2^2}{4a_4}\right)-\left(\sqrt{a_4} r^2 + \frac{a_2}{\sqrt{4a_4}}\right)^2
\right]^{3/2} rdr
\end{eqnarray}

{It is useful to calculate the integral
\begin{equation}
J(x) = \int (\lambda^2 - x^2)\, dx, \quad \lambda >0.
\end{equation}
After elementary integration by parts, we obtain:
\begin{equation}
J(x) =\frac{3\lambda^4}{8} \arcsin\left(\frac{x}{\lambda}\right)
+ \frac{3\lambda^2}{8} x \left(\lambda^2 - x^2\right)^{1/2}
+ \frac{1}{4} x \left(\lambda^2 - x^2\right)^{3/2},
\end{equation}
or in the more useful form:
\begin{equation}
\label{eq-app-Jx}
J(x) =\lambda^4\left[\frac{3}{8} \arcsin\left(\frac{x}{\lambda}\right)
+ \frac{3}{8} \left(\frac{x}{\lambda}\right) \left(1 - \left(\frac{x}{\lambda}\right)^2\right)^{1/2}
+ \frac{1}{4} \left(\frac{x}{\lambda}\right) \left(1 - \left(\frac{x}{\lambda}\right)^2\right)^{3/2}\right].
\end{equation}
We also notice that:
\begin{equation}
\label{eq-app-Jx-def}
J(\lambda)=\frac{3\pi}{16} \lambda^4, \quad J(0)=0, \quad J(-\lambda)= -\frac{3\pi}{16} \lambda^4.
\end{equation}
}
Using now the notation
\begin{equation}
\lambda = \sqrt{\rho_{0}+\frac{a_2^2}{4a_4}},
\end{equation}
and the change of variables
\begin{eqnarray} \vspace{0.2cm}
\displaystyle u = \sqrt{a_4} r^2 + \frac{a_2}{\sqrt{4a_4}}, \quad du = 2 \sqrt{a_4} dr,\\ \vspace{0.2cm}
\displaystyle r=0 \Longrightarrow u_0=\frac{a_2}{\sqrt{4a_4}}, \\
\displaystyle r=R_\perp \Longrightarrow u_\perp=\lambda,
\end{eqnarray}
our integral becomes:
\begin{equation}
I=\frac{4\pi}{3\sqrt{a_z a_4}} \int_{u_0}^{\lambda} (\lambda^2-u^2)^{3/2} du =
\frac{4\pi}{3\sqrt{a_z a_4}}\left(J(\lambda)-J(u_0)\right).
\end{equation}
Introducing the parameter:
\begin{equation}
\label{eq-app-q3D-eta}
\eta=\frac{a_2}{\sqrt{4a_4 \rho_{0}}}>0, \Longrightarrow \frac{u_0}{\lambda}=\frac{\eta}{\sqrt{1+\eta^2}},
\end{equation}
and using \eqref{eq-app-Jx} and \eqref{eq-app-Jx-def}, we finally obtain:
\begin{eqnarray} \vspace{0.2cm} \nonumber
I&=&\displaystyle \frac{4\pi}{3\sqrt{a_z a_4}} \lambda^4
\left[\frac{3}{8}\left(\frac{\pi}{2}-\arcsin\frac{\eta}{\sqrt{1+\eta^2}}
\right)
-\frac{3}{8} \frac{\eta}{\sqrt{1+\eta^2}} \frac{1}{(1+\eta^2)^{1/2}}
-\frac{1}{4} \frac{\eta}{\sqrt{1+\eta^2}} \frac{1}{(1+\eta^2)^{3/2}}
\right],\\
&=& \displaystyle \frac{4\pi}{3\sqrt{a_z a_4}} \lambda^4
\left[\frac{3}{8}\arccos\frac{\eta}{\sqrt{1+\eta^2}}
-\frac{3}{8} \frac{\eta}{1+\eta^2}
-\frac{1}{4} \frac{\eta}{(1+\eta^2)^2}
\right].
\end{eqnarray}
Using that:
\begin{equation}
\lambda^4 = \rho_{0}^2 (1+\eta^2)^2 = \frac{a_2^4}{(4a_4)^2} \frac{(1+\eta^2)^2}{\eta^4},
\end{equation}
we obtain a non-linear equation in $\eta$:
\begin{eqnarray} \vspace{0.2cm} \nonumber
I =C_\text{S} &=& \frac{8\pi a_2^4}{3 a_z^{1/2}\, (4 a_4)^{5/2}} \frac{1}{\eta^4}
\left[\frac{3(1+\eta^2)^2}{8}\arccos\frac{\eta}{\sqrt{1+\eta^2}}
-\frac{3}{8} \eta (1+\eta^2)
-\frac{1}{4} \eta
\right],\\
&=& \frac{\pi a_2^4}{a_z^{1/2}\, (4 a_4)^{5/2}} \frac{1}{\eta^4}
\left[(1+\eta^2)^2 \arccos\frac{\eta}{\sqrt{1+\eta^2}}
-\eta^3 -\frac{5}{3} \eta
\right].
\end{eqnarray}
\medskip

To summarize this case, we have to
\begin{itemize}
  \item find the root $\eta >0$ of the non-linear equation:
  \begin{equation}
  \label{eq-app-q3D-feta}
  f(\eta)=\underbrace{\frac{a_z^{1/2}\, (4 a_4)^{5/2}}{\pi a_2^4} C_\text{S}}_{A_\eta}\, \eta^4
  -(1+\eta^2)^2 \arccos\frac{\eta}{\sqrt{1+\eta^2}}
  +\eta^3 +\frac{5}{3} \eta =0,
  \end{equation}
  \begin{equation}
  \label{eq-app-q3D-fetap}
  f'(\eta)=4A_\eta \eta^3
  -4\eta(1+\eta^2) \arccos\frac{\eta}{\sqrt{1+\eta^2}}
  +(1+\eta^2)+3\eta^2 +\frac{5}{3} ,
  \end{equation}

  \item compute
  \begin{equation}
  \label{eq-app-q3D-rhof}
  \rho_{0}=\frac{a_2^2}{4a_4 \eta^2},
  \end{equation}
  and the dimensions of the condensate:
  \begin{equation}
  \label{eq-app-q3D-rperp}
  R_\perp^2= \frac{-a_2 + \sqrt{a_2^2+4\rho_{0} a_4}}{2 a_4},
  \end{equation}
  \begin{equation}
  \label{eq-app-q3D-zmax}
  R_{zmax} = z|_{r=0}=\left(\frac{\rho_{0}}{a_z}\right)^{1/2}.
  \end{equation}
\end{itemize}

\subsubsection{Case $a_2=0$: pure quartic potential}
The integration is carried exactly in the same manner, the difference coming from the limits of the integration following $r$. We obtain
 \begin{equation}
I=
\frac{4\pi}{3\sqrt{a_z a_4}}\left(J(\lambda)-J(0)\right)
=\frac{4\pi}{3\sqrt{a_z a_4}} \, \frac{3\pi}{16} \lambda^4,
\end{equation}
with $\lambda=\sqrt{\rho_{0}}$. Finally
\begin{equation}
I =C_\text{S} = \frac{\pi^2}{2 \sqrt{a_z 4a_4}} \left(\rho_{0}\right)^2,
\end{equation}
and
\begin{equation}
\rho_{0}= \frac{1}{\pi}\left(2 a_z^{1/2}\, (4a_4)^{1/2}\right)^{1/2} \,C_\text{S}^{1/2}.
\end{equation}
\begin{equation}
R_\perp = \left(\frac{\rho_{0}}{a_4}\right)^{1/4}, \quad R_{zmax} = \left(\frac{\rho_{0}}{a_z}\right)^{1/2}.
\end{equation}
\subsubsection{Case $a_2<0$: quartic-quadratic potential}

For this case, $a_2<0, a_z>0, a_4>0$ and we distinguish two subcases:
\begin{itemize}

\item If $\rho_{0}<0$, the condensate has a hole.
We infer from \eqref{eq-app-q3D-roots} that there are two roots  $R_\perp^\pm$
\begin{equation}
(R^\pm_\perp)^2 = \frac{-a_2 \pm \sqrt{a_2^2+4\rho_{0} a_4}}{2 a_4} > 0.
\end{equation}
and the condensate extends in central plane from $R_\perp^-$ to $R_\perp^+$.

The integration is carried exactly in the same manner, the difference coming from the limits of the integration following $r$. We obtain
 \begin{equation}
I=
\frac{4\pi}{3\sqrt{a_z a_4}}\left(J(\lambda)-J(-\lambda)\right)
=\frac{4\pi}{3\sqrt{a_z a_4}} \, 2\frac{3\pi}{16} \lambda^4.
\end{equation}
and finally
\begin{equation}
I =C_\text{S} = \frac{\pi^2}{\sqrt{a_z 4a_4}} \left(\rho_{0}+\frac{a_2^2}{4a_4}\right)^2.
\end{equation}
The value of $\rho_{0}$ results as:
\begin{equation}
\rho_{0}= \frac{a_z^{1/4}\, (4a_4)^{1/4}}{\pi} \,C_\text{S}^{1/2} - \frac{a_2^2}{4a_4} = \frac{a_2^2}{4a_4} (\xi-1).
\end{equation}
Since $\rho_{0}<0$, this case is obtained if:
\begin{equation}
\xi = \frac{a_z^{1/4}\, (4a_4)^{5/4}}{\pi a_2^2} \,(C_\text{S})^{1/2} = \frac{\sqrt{A_\eta}}{\sqrt{\pi}} < 1.
\end{equation}
The dimensions of the condensate are
  \begin{equation}
  R_{max}^2= \frac{-a_2 + \sqrt{a_2^2+4\rho_{0} a_4}}{2 a_4},
  \end{equation}
  \begin{equation}
 R_{zmax} =z{\big |}_{r^2=\frac{-a_2}{2a_4}}=\frac{1}{\sqrt{a_z}}
  \sqrt{\rho_{0}+\frac{a_2^2}{4a_4}}.
  \end{equation}

\item If $\rho_{0}>0$, the condensate has only a depletion centered around $z=0$ (the density profile has not any more the maximum at $z=0$). This case occurs when:
    \begin{equation}
    \xi = \frac{a_z^{1/4}\, (4a_4)^{5/4}}{\pi a_2^2} \,C_\text{S}^{1/2} = \frac{\sqrt{A_\eta}}{\sqrt{\pi}} > 1.
    \end{equation}
    The computation is the same as for the case of "quartic + quadratic" potential, with the difference that the root $\eta$ is now negative. In particular
      \begin{equation}
 R_{zmax} =\frac{1}{\sqrt{a_z}}
  \sqrt{\rho_{0}+\frac{a_2^2}{4a_4}}.
  \end{equation}
\end{itemize}

\subsubsection{Summary for the 3D combined quartic and quadratic potential}
\vspace{-0.5cm}
{\small
\[
\vadim{V}^\text{eff} = \frac12\left(a_2 r^2 + a_4 r^4 + a_z z^2\right).
\]
\vspace{-1cm}
\[
\left[
\begin{array}{l}   
\text{Compute} \quad C_\text{S} = 2\varepsilon^2C_g,\\
\mbox{$\bullet$ if $a_2 =0$, } \\ \vspace{0.15cm}
\left\{\begin{array}{l}
\displaystyle \rho_{0} = \frac{1}{\pi}\left(2 a_z^{1/2}\, (4a_4)^{1/2}\right)^{1/2} \,C_\text{S}^{1/2},\\ \vspace{0.15cm}
\displaystyle   R_{max}= \left(\frac{\rho_{0}}{a_4}\right)^{1/4}, \quad
\displaystyle R_{zmax} = \left(\frac{\rho_{0}}{a_z}\right)^{1/2}.
\end{array}
\right. \\ \vspace{0.15cm}
\mbox{$\bullet$ else} \\ \vspace{0.15cm}
\begin{array}{l} 
--> \mbox{compute} \\ \vspace{0.15cm}
\displaystyle A_{\eta} = \frac{a_z^{1/2}\, (4 a_4)^{5/2}}{\pi a_2^4} C_\text{S}\\ \vspace{0.15cm}
--> \mbox{define the function} \\ \vspace{0.15cm}
f(\eta)={A_\eta}\, \eta^4
  -(1+\eta^2)^2 \arccos\frac{\eta}{\sqrt{1+\eta^2}}
  +\eta^3 +\frac{5}{3} \eta \\ \vspace{0.15cm}
\left[ 
\begin{array}{l}\vspace{0.15cm}
\mbox{$\bullet$ if $a_2 >0$, } \\ \vspace{0.15cm}
---->\mbox{find the positive root $\eta \in [0,200]$ of $f(\eta)=0$} \\ \vspace{0.15cm}
\left\{\begin{array}{l}
\displaystyle \rho_{0} =\frac{a_2^2}{4a_4 \eta^2},\\ \vspace{0.15cm}
\displaystyle   R_{max}= \left(\frac{-a_2 + \sqrt{a_2^2+4\rho_{0} a_4}}{2 a_4}\right)^{1/2}, \quad
\displaystyle R_{zmax} = \left(\frac{\rho_{0}}{a_z}\right)^{1/2}.
\end{array}
\right. \\ \vspace{0.15cm}
\mbox{$\bullet$ else } \\ \vspace{0.15cm}
\displaystyle \xi = \frac{\sqrt{A_\eta}}{\sqrt{\pi}}, \\ \vspace{0.15cm}
\left[\begin{array}{l} 
\mbox{$\bullet$ if $\xi < 1$, } \\ \vspace{0.15cm}
\left\{\begin{array}{l}
\displaystyle \rho_{0} =\frac{a_2^2}{4a_4} (\xi-1),\\ \vspace{0.15cm}
\displaystyle   R_{max}= \left(\frac{-a_2 + \sqrt{a_2^2+4\rho_{0} a_4}}{2 a_4}\right)^{1/2}, \quad
\displaystyle R_{zmax} = \left(\frac{\rho_{0}}{a_z}+\frac{a_2^2}{4a_4 a_z}\right)^{1/2}.
\end{array}
\right. \\ \vspace{0.15cm}
\mbox{$\bullet$ else} \\ \vspace{0.15cm}
\mbox{find the {\bf negative} root $\eta \in [-200, 0]$ of $f(\eta)=0$} \\ \vspace{0.15cm}
\left\{\begin{array}{l}
\displaystyle \rho_{0} =\frac{a_2^2}{4a_4 \eta^2},\\ \vspace{0.15cm}
\displaystyle   R_{max}= \left(\frac{-a_2 + \sqrt{a_2^2+4\rho_{0} a_4}}{2 a_4}\right)^{1/2}, \quad
\displaystyle R_{zmax} = \left(\frac{\rho_{0}}{a_z}+\frac{a_2^2}{4a_4 a_z}\right)^{1/2}.
\end{array}
\right.
\end{array}\right. 
\end{array}
\right. 
\end{array}
\end{array}
\right. 
\]
}

\bibliographystyle{elsarticle-num}

\end{document}